\documentclass[11pt]{article}

\usepackage{amsthm,amsfonts,amssymb,amsmath}
\usepackage{enumerate}
\usepackage{fullpage}
\usepackage{tikz-cd}
\usepackage{indentfirst}
\usepackage[colorlinks=true]{hyperref}
\newtheorem{thm}{Theorem}[section]

\newtheorem{theorem}[thm]{Theorem}
\newtheorem{proposition}[thm]{Proposition}
\newtheorem{lemma}[thm]{Lemma}
\newtheorem{corollary}[thm]{Corollary}

\newtheorem{conjecture}[thm]{Conjecture}

\theoremstyle{definition}
\newtheorem{definition}[thm]{Definition}

\theoremstyle{remark}
\newtheorem{remark}[thm]{Remark}

\newcommand{\Span}{{\mathrm{Span}}}
\newcommand{\trdeg}{{\mathrm {trdeg}}}

\newcommand{\Alb}{{\mathrm {Alb}}}
\newcommand{\Bs}{{\mathrm{Bs}}}

\newcommand{\len}{{\mathrm{length}}}

\newcommand{\SBs}{{\mathrm {SBs}}}
\newcommand{\ddiv}{{\mathrm {div}}}

\newcommand{\Coker}{{\mathrm{Coker}}}

\renewcommand{\div}{{\mathrm{div}}}

\newcommand{\ch}{{\mathrm{CH}}}

\newcommand{\NS}{{\mathrm{NS}}}

\newcommand{\Pic}{\mathrm{Pic}}

\newcommand{\supp}{{\mathrm{Supp}}}

\newcommand{\vol}{{\mathrm{vol}}}

\theoremstyle{plain}
\newtheorem{introtheorem}{Theorem}

\DeclareMathOperator{\Spec}{\textnormal{Spec}}

\DeclareMathOperator{\ord}{\textnormal{ord}}
\DeclareMathOperator{\Gal}{\textnormal{Gal}}

\begin{document}

\title{Arithmetic degrees of dynamical systems over fields of positive characteristic}
\author{Yubo Zhou}


\maketitle






\begin{abstract} 
In this article, we study  heights and arithmetic degrees on dynamical systems $(X,f)$ over an algebraically closed field $\mathbf{k}$ of positive characteristic. By taking
spread-outs, we may choose a finitely generated field such that the system can be defined.
Over such fields, we give a direct proof of the Northcott property for Moriwaki heights
in this setting. These
heights allow us to define the arithmetic degrees
of $x$ independently of the choice of the spread-out. We also define the arithmetic dynamical degrees of $f$ on its projective models.
 The intersection inequalities on these
models yield the formula
$\widehat{\lambda}_i(f)=\max\{\lambda_i(f),\lambda_{i-1}(f)\}$,
from which we deduce the fundamental inequality
$\overline{\alpha}_f(x)\leq\lambda_1(f)$ for every $x\in X_f(\mathbf{k})$. We then formulate the Kawaguchi--Silverman
conjecture in this setting. When $f$ is a surjective endomorphism
and $\mathbf{k}$ has infinite transcendence degree,
we prove the predicted equality $\alpha_f(x)=\lambda_1(f)$ for every
 purely transcendental point.
\end{abstract}

\tableofcontents

\section{Introduction}

Let $\mathbf{k}$ be a algebraically closed field with prime field $k$. Let $X$ be a projective variety over $\mathbf{k}$, and let
$f:X\dashrightarrow X$ be a dominant rational self-map. By spreading out, we may choose a subfield $K\subseteq\mathbf{k}$, finitely generated over $k$,
over which both $X$ and $f$ can be defined. 

Passing to a finitely generated field $K$ allows us to study the arithmetic of $(X,f)$
through its models. Over number fields, Song \cite{Song2026}
established a
higher-codimensional  Yuan's inequality and used it to defined arithmetic dynamical degree.
Luo and Song \cite{LS2025} extended this formula to finitely
generated extensions of $\mathbb{Q}$. Using Moriwaki heights
over such fields of definition, they also defined the upper and
lower arithmetic degrees of points for dynamical systems over
arbitrary fields of characteristic zero and proved that these
degrees are independent of the chosen field of definition.

In this paper, we study both the arithmetic degrees of points
and the arithmetic dynamical degrees of maps in positive
characteristic. We begin with the height theory over finitely
generated extensions of $\mathbb{F}_p$. We give a direct proof of the Northcott property for the heights considered here.

\begin{introtheorem}
Let $K/\mathbb{F}_q$ be a finitely generated extension
of transcendence degree $d\geq1$, and let
$\mathcal{B}=(B;H_1,\dots,H_{d-1})$ be a big polarization.
Let $X$ be a projective variety over $K$, and let $L$
be an ample line bundle on $X$.
Then, for every real number $M$ and every positive
integer $e$, the set
\begin{equation*}
    \left\{
        P\in X(\overline{K})
        \mid h_L^{\mathcal{B}}(P)\leq M,\quad
        [K(P):K]\leq e
    \right\}
\end{equation*}
is finite.
\end{introtheorem}
We use these heights to measure the growth of the height of a point
along its orbit. With the notation of the preceding theorem, assume
that $X$ is normal, and let $f:X\dashrightarrow X$ be a dominant
rational self-map over $K$. Denote by $X_f(\overline{K})$ the set
of points in $X(\overline{K})$ whose forward orbits under $f$ are
well defined. 
For $x\in X_f(\overline{K})$, one defined the arithmetic
degrees of $x$ under $f$ are defined by
\begin{equation*}
     \overline{\alpha}_f(x)
    :=\limsup_{m\to\infty} h^+(f^m(x))^{1/m}
\end{equation*}
\begin{equation*}
    \underline{\alpha}_f(x)
    :=\liminf_{m\to\infty} h^+(f^m(x))^{1/m}.
\end{equation*}
When these two limits agree, we denote their common value by
$\alpha_f(x)$ and call it the arithmetic degree of $x$ under $f$. We will show that the  arithmetic degree only  depends the dynamical system.

Another important invariant of $f$ is the dynamical degree
.Let $f:X\dashrightarrow X$ be a dominant rational self-map of a normal projective variety $X$ over $\mathbf{k}$.
Consider the diagram
\begin{equation*}
    \begin{tikzcd}
    & \Gamma_{f^m}
      \arrow[dl,"\pi_{m,1}"']
      \arrow[dr,"\pi_{m,2}"] & \\
    X \arrow[rr,dashed,"f^m"] && X
\end{tikzcd}
\end{equation*}
Where $\Gamma_{f^m}$ is the normalization of the closure of the
graph of $f^m$ in $X\times X$. The first projection $\pi_{m,1}$
is birational, and the second projection $\pi_{m,2}$ is dominant
and generically finite. We fixed a big and nef line bundle $L$ on $X$, one can define the $i$-th dynamical degree of $f$ as
\begin{equation*}
        \lambda_i(f)
    =
    \lim_{m\to\infty}
    \left((f^m)^*L^i\cdot L^{\dim X-i}\right)^{1/m},
    \qquad 0\leq i\leq \dim X.
\end{equation*}
Dinh and Sibony \cite{DS05} proved the existence of this
limit in the complex
 setting, and Truong \cite{Tru20}
gave an algebraic proof in arbitrary characteristic.
In our setting, Dang \cite{Dang20} showed that the degrees
$\deg_{i,L}(f^m)$ are submultiplicative with a constant independent
of $m$, so the limit therefore exists by Fekete's lemma. 

Motivated by the work of \cite{Song2026} and 
\cite{LS2025}, we study arithmetic dynamical degrees over finitely
generated fields of positive characteristic. In this setting, a
projective model is a variety over a finite field, and the arithmetic
degrees of the iterates are given by intersection numbers on the
model. To estimate these intersection numbers, we use the following
algebraic form of Siu's inequality.

\begin{introtheorem}
    Let $\mathcal{X}\to k$ be a normal projective variety with dimension $n$, and $0\leq r\leq n-1$  be an integer. Let $\mathcal{L}$, $\mathcal{M}$, $\mathcal{N}_i$, $i=1,2,\dots, r$ be nef line bundles such that 
      \begin{equation*}
      \left(\mathcal{M}^{n-r}\cdot \mathcal{N}_1\cdots \mathcal{N}_r\right)>\left(n-r\right) \left(\mathcal{M}^{n-r-1}\cdot \mathcal{L}\cdot \mathcal{N}_1\cdots \mathcal{N}_r\right).
  \end{equation*}
Then we have 
    \begin{equation*}
      \mathcal{L}\cdot \mathcal{N}_1\cdots \mathcal{N}_r\leq \mathcal{M}\cdot \mathcal{N}_1\cdots \mathcal{N}_r.
  \end{equation*}
\end{introtheorem}

We apply this inequality to the finitely generated field 
$K\subseteq\mathbf{k}$ chosen above for $(X,f)$. Let $K/\mathbb{F}_q$ be a finitely
generated field of transcendence degree $d\geq 1$, and let
$f:X\dashrightarrow X$ be a dominant rational self-map of a normal
projective variety of dimension $n$ over $K$. Choose a normal
projective model $B$ of $K$ and a normal projective model
$\pi:\mathcal{X}\to B$ of $X$. Then $f$ induces a rational
self-map $F:\mathcal{X}\dashrightarrow\mathcal{X}$ over $B$.
\begin{center}
   \begin{tikzcd}
 & \Gamma_{F} \arrow{dl}[swap]{\pi_{1}} \arrow{dr}{\pi_2}\\
\mathcal{X} \arrow[dashed]{rr}{F} && \mathcal{X} ,
\end{tikzcd} 
\end{center} 
For $1\leq i\leq n$, the preceding intersection inequality gives
a submultiplicative estimate for the $i$-th degrees of the iterates
$F^m$. 
So we can use it to defines the $i$-th arithmetic dynamical
degree $\widehat{\lambda}_i(f)$. The following theorem expresses
$\widehat{\lambda}_i(f)$ in terms of the dynamical degrees on the
generic fiber.

\begin{introtheorem}
 Let $f:X\dashrightarrow X$ be a dominant rational self-map,
and let $1\leq k\leq n$ be an integer.
With the notation defined as above, we have
  \begin{equation*}
  \widehat{\lambda}_k\left(f\right)=\max\{\lambda_k(f), \lambda_{k-1}(f)\}.
  \end{equation*}
\end{introtheorem}
Moreover, over number field, the Kawaguchi--Silverman conjecture
say that the arithmetic degree $\alpha_f(x)$ exists for
$x\in X_f(\overline{\mathbb{Q}})$ and that
$\alpha_f(x)=\lambda_1(f)$ when the orbit is Zariski dense
\cite{KS16a}. Our aim is to study this question
over algebraically closed fields of positive characteristic.

Following the approach of Luo and Song
\cite{LS2025} in characteristic zero, we choose a finitely generated field 
 $K$ which the system $(X,f)$ can be defined. Comparing the growth of heights along the orbit of $x$
with the first arithmetic dynamical degree of $f$, we obtain
the fundamental inequality.
\begin{introtheorem}
For every $x\in X_f(\mathbf{k})$, we have
\begin{equation*}
    \overline{\alpha}_f(x)\leq\lambda_1(f).
\end{equation*}
\end{introtheorem}
A direct analogue of the classical Kawaguchi--Silverman conjecture
in positive characteristic would take the following form.

\begin{conjecture}
Let $X$ be a normal projective variety over a field
$\mathbf{k}$ 
and  $f:X\dashrightarrow X$ be a dominant
rational self-map. Then 
\begin{itemize}
    \item
    The arithmetic degree $\alpha_f(x)$ exists
    for every $x\in X_f(\mathbf{k})$.

    \item
    For any $x\in X_f(\mathbf{k})$, if the $f$-orbit $O_f(x)$ is Zariski dense,
   then $
        \alpha_f(x)=\lambda_1(f).$
\end{itemize}
\end{conjecture}
When $f$ is an endomorphism, we prove that $\alpha_f(x)$ exists
for every $x\in X(\mathbf{k})$. For surjective endomorphisms,
the following theorem gives the equality $\alpha_f(x)=\lambda_1(f)$
at purely transcendental points.

\begin{introtheorem}
Let $X$ be a normal projective variety over $\mathbf{k}$ and let $f:X\to X$ be a surjective
endomorphism.
Suppose that $(X,f)$ can be defined over a finitely generated field $K\subseteq\mathbf{k}$.
Then, for any purely transcendental point
$x\in X(\mathbf{k})$ with respect to
$X_{\overline{K}}\to\Spec(\overline{K})$, we have
$\alpha_f(x)=\lambda_1(f)$.
In particular, there exists a countable collection $S$
of proper closed subvarieties of $X$ such that
$\alpha_f(x)=\lambda_1(f)$
for every
$x\in X(\mathbf{k})\smallsetminus
\bigcup_{V\in S}V(\mathbf{k})$.
\end{introtheorem}

\subsection*{Organization of the paper}
In Section \ref{height chatper}, we study Moriwaki heights over finitely generated fields of
positive characteristic. We give a direct proof of the Northcott
property for big polarizations and compare the heights under the field extension.

In Section~\ref{fg case chatper}, we use intersection inequality 
to prove the existence of arithmetic dynamical degrees $\widehat{\lambda}_1(f)$. We then prove the prove the relative degree formula and establish the fundamental inequality $\overline{\alpha}_f(x)\leq\lambda_1(f)$  over
finitely generated fields.

In Section~\ref{algebraic closed case}, we define the arithmetic degrees of points over
algebraically closed fields of positive characteristic by choosing a finitely generated field. 
We prove it's  independent of the choose of the field. Moreover, we study the Kawaguchi--Silverman conjecture in this
setting and conclude with our result on purely transcendental
points for surjective endomorphisms.
\subsection*{Acknowledgments}
This work forms part of my master's thesis at Harbin Institute of
Technology. I would like to express my great gratitude to my advisor Chunhui Liu for suggesting this problem
and for guiding me throughout this project. I also thank Wenbin Luo
and Jiarui Song for introducing me to their work  at the 2024 workshop on
Arakelov geometry and Nevanlinna theory in Harbin.
 In addition, I thank Heer Zhao for funding my access to ChatGPT
while I was making the final revisions to this paper.

\section{Heights over finitely generated fields of positive characteristic}\label{height chatper}

In this section, we study Moriwaki heights over finitely
generated fields of positive characteristic.
We establish the Northcott property for big polarizations
and compare heights under finitely generated field extensions.
\subsection{Moriwaki height and the Northcott property}
Let $K$ be a finitely generated field over
$k=\mathbb{F}_q$ with $\trdeg_k(K)=d\geq1$,
and fix an algebraic closure $\overline{K}$.
Choose a polarization
$\mathcal{B}=(B;H_1,\dots,H_{d-1})$,
where $B$ is a normal projective $k$-variety with
function field $K$, and the $H_i$ are nef line bundles.

For $P=[x_0:\cdots:x_n]\in\mathbb{P}^n(K)$,
the \emph{naive height} is defined by
\begin{equation*}
    h^{\mathcal{B}}_{\mathrm{nv},K}(P)
    :=
    \sum_{\Gamma}
    \max_i\{-\ord_{\Gamma}(x_i)\}
    \left(\Gamma\cdot H_1\cdots H_{d-1}\right),
\end{equation*}
where $\Gamma$ runs over the prime divisors of $B$
.
The product formula shows that this is independent
of the homogeneous coordinates.

For a finite extension $K'/K$, let $\nu:B'\to B$
be the normalization of $B$ in $K'$ and put
$\mathcal{B}'=(B';\nu^*H_1,\dots,\nu^*H_{d-1})$.
For $P\in\mathbb{P}^n(K')$, set
\begin{equation*}
    h^{\mathcal{B}}_{\mathrm{nv},K}(P)
    :=
    \frac{1}{[K':K]}
    h^{\mathcal{B}'}_{\mathrm{nv},K'}(P).
\end{equation*}
By the projection formula, this defines a height on
$\mathbb{P}^n(\overline{K})$ independent of $K'$.

More generally, let $X$ be a  projective
$K$-variety and $L$ a line bundle on $X$.
Choose a $B$-model $(\mathcal{X},\mathcal{L})$ of $(X,L)$:
a integral projective scheme
$f:\mathcal{X}\to B$ with a line bundle $\mathcal{L}$
such that $\mathcal{X}_K\simeq X$ and
$\mathcal{L}_K\simeq L$. In this case, the Moriwaki height $h^{\mathcal{B}}_{(\mathcal{X},\mathcal{L})}$
is given by the intersection theory.

\begin{definition}
\label{defi of moriwaki height}
For $P\in X(\overline{K})$, let $\Delta_P$ be the
Zariski closure of its image in $\mathcal{X}$. 
The \emph{Moriwaki height} associated with
$\mathcal{B}$ and $(\mathcal{X},\mathcal{L})$ is
\begin{equation*}
    h^{\mathcal{B}}_{(\mathcal{X},\mathcal{L})}(P)
    :=
    \frac{
       \deg \left(
            \mathcal{L}|_{\Delta_P}
            \cdot f^*H_1|_{\Delta_P}
            \cdots f^*H_{d-1}|_{\Delta_P}
        \right)
    }{
        [K(P):K]
    }.
\end{equation*}
\end{definition}

We next show that, for fixed $\mathcal{B}$ and $(X,L)$,
the height is independent of the model up to $O(1)$.
We first record the compatibility of heights with finite
extensions.

\begin{proposition}\label{height under finite extension}
Let $K'/K$ be a finite extension, and let $\pi:B'\to B$
be the normalization of $B$ in $K'$.
Put $\mathcal{B}'
    =(B';\pi^*H_1,\dots,\pi^*H_{d-1}),
     $. 
    Let $\mathcal{X}'$ be a main component
of $X_{K'}$ in $\mathcal{X}\times_B B'$ 
whose generic fiber
contains the point induced by $P$.and denote the projection by
$\pi':\mathcal{X}'\to\mathcal{X}$.
Then
\begin{equation*}
    h^{\mathcal{B}'}_{(\mathcal{X}',\pi'^*\mathcal{L})}(P)
    =
    [K':K]\,
    h^{\mathcal{B}}_{(\mathcal{X},\mathcal{L})}(P).
\end{equation*}
for every $P\in X(\overline{K})$.

\end{proposition}

\begin{proof}
For any $P\in X(\overline{K})$. Let $\Delta_P$(resp.$\Delta_{P}^{'}$) be the closure in $\mathcal{X}$(resp.$\mathcal{X}^{'}$), then we have
  \begin{equation*}
        \begin{array}{ll}
            h^{\left(B';\pi^{*}\left(H_1\right),\dots,\pi^{*}\left(H_{d-1}\right)\right)}_{(\mathcal{X}',\pi^{\prime *}\left(\mathcal{L}\right))}(P)&\\
            &=\dfrac{\deg\left(\left(\pi^{\prime *}f^{*}H_1|_{\Delta_P^{\prime}}\right)\cdots \left(\pi^{\prime *}f^{*}H_{d-1}|_{\Delta_P^{\prime}}\right)\cdot \left(\pi^{\prime *}\mathcal{L}|_{\Delta_P^{\prime}}\right) \right)}{\deg\left(\Delta_P^{'}\to B'\right)}\\
            &=\dfrac{\deg\left(\Delta_P^{'}\to \Delta_P\right)}{\deg\left(\Delta_P^{'}\to B'\right)}\left(\left(f^{*}H_1|_{\Delta_P}\right)\cdots \left(f^{*}H_{d-1}|_{\Delta_P}\right)\cdot \left(\mathcal{L}|_{\Delta_P}\right) \right)\\
            &=[K':K]h^{\left(B;H_1,\dots,H_{d-1}\right)}_{(\mathcal{X},\mathcal{L})}(P).
        \end{array}
    \end{equation*}
    ~
\end{proof}

The naive height is recovered from the standard model
of projective space.

\begin{proposition}
\label{naive height given by moriwaki height}
On $\mathbb{P}^n(\overline{K})$, we have 
$h_{\mathrm{nv}}^{\mathcal{B}}(\cdot)
    =
    h^{\mathcal{B}}_{(\mathbb{P}^n_B,\mathcal{O}(1))}(\cdot).$

\end{proposition}

\begin{proof}
  By proposition \ref{height under finite extension}, it remains to prove for  $P\in X(K)$.
Fix a basis $\{X_0,\dots,X_n\}$ of $H^0\left(\mathbb{P}^n_k,\mathcal{O}(1)\right)$. Let $\Delta_P$ be the closure of $P$ in $\mathbb{P}^n_B$.
Replacing $B$ by the normalization of $\Delta_P$, we may assume that $P$ extends to a
section  $s_P:B\to  \mathbb{P}^n_B$. We may assume $s^{*}_P(X_0)\neq 0$. If we set $a_i=s^{*}_P(X_i)/s^{*}_P(X_0)$ for $i=0,\dots,n$, then $a_i\in K$ and $h_{nv}^{\mathcal{B}}(P)=h_{nv}^{\mathcal{B}}\left([a_0:a_1:\cdots:a_n]\right)$. Since $\ddiv\left(s_P^{*}(X_0)\right)=\sum_{\Gamma}\ord_{\Gamma}\left(s^{*}_P(X_0)\right)\Gamma$, we have
  \begin{equation*}
        \begin{array}{ll}
       \deg\left(\left(f^{*}H_1\right)\cdots
        \left(f^{*}H_{d-1}\right)\right.
        &\cdot\left.\mathcal{O}(1)\cdot\Delta_P\right)\\
        &=\deg\left(
            \left(s_P^{*}f^{*}H_1\right)\cdots
            \left(s_P^{*}f^{*}H_{d-1}\right)
            \cdot\left(s_P^{*}\mathcal{O}(1)\right)
            \cap[B]
        \right)\\
        &=\deg\left(
            H_1\cdots H_{d-1}\cdot s_P^{*}\left(\mathcal{O}(1)\right)\cap[B]
        \right)\\
        &=\sum\limits_{\Gamma}\ord_{\Gamma}(s_P^{*}(X_0))
        \deg\left(H_1|_{\Gamma}\cdots H_{d-1}|_{\Gamma}\right).
        \end{array}
    \end{equation*}
On the other hand, since $s_P^{*}(X_i)$ generate $s_P^{*}\left(\mathcal{O}(1)\right)$. We have
\begin{equation*}
    \begin{array}{ll}
        \ord_{\Gamma}(s_P^{*}(X_0))&=\len_{\Gamma}\left(\frac{s_P^{*}(\mathcal{O}(1))}{s_P^{*}(X_0)}\right)\\
        &=\len_{\Gamma}\left(\frac{\mathcal{O}_B s_P^{*}(X_0)+\cdots +\mathcal{O}_B s_P^{*}(X_n)}{s_P^{*}(X_0)}\right)\\
        &=\len_{\Gamma}\left(\frac{\mathcal{O}_B a_0+\cdots +\mathcal{O}_B a_n}{\mathcal{O}_B}\right)\\
              &=\max\limits_i\left\{-\ord_{\Gamma}(a_i)\right\}.
    \end{array}
\end{equation*}
\end{proof}

\begin{proposition}\label{upper boundary of height}
Let $\Bs(L)$ be the base locus of $L$.
Then function $h^{\mathcal{B}}_{(\mathcal{X},\mathcal{L})}(\cdot)$
is bounded below on
$(X\smallsetminus\Bs(L))(\overline{K})$.
\end{proposition}

\begin{proof}
Choose an ample line bundle $A$ on $B$ such that
$f_*\mathcal{L}\otimes A$ is generated by global sections. So the map $ H^0(\mathcal{X},\mathcal{L}\otimes f^{*}A)\otimes \mathcal{O}_B \to f_{*}\mathcal{L}\otimes A$  is surjective. Hence, $ H^0(\mathcal{X},\mathcal{L}\otimes f^{*}A)\otimes \mathcal{O}_{\mathcal{X}} \to f^{*}(f_{*}\mathcal{L})\otimes f^{*}A$ is also surjective by pullback. Set
\begin{equation*}
    \mathcal{Z}:=\supp\left(\Coker\left(H^0(\mathcal{X},\mathcal{L}\otimes f^{*}A)\otimes \mathcal{O}_{\mathcal{X}}\to \mathcal{L}\otimes f^{*}A \right)\right),
\end{equation*}
then we have $\mathcal{Z}_K=\Bs(L)$. Take a basis $\{s_0,\dots,s_k\}$ of $H^0(\mathcal{X},\mathcal{L}\otimes f^{*}A)$, for any $P\in(X\smallsetminus\Bs(L))(\overline{K})$. There exists a $s_i$ with $s_i(P)\neq 0$, which shows that $\mathcal{L}\otimes f^{*}A|_{\Delta_P}$ is effective. Thus
\begin{equation*}
    \deg\left((f^{*}H_1|_{\Delta_P})\cdots (f^{*}H_{d-1}|_{\Delta_P})\cdot (\mathcal{L}\otimes f^{*}A|_{\Delta_P})\right)\geqslant 0.
\end{equation*}
By projection formula, we get
\begin{equation*}
  \begin{array}{rl}
       \deg  \left((f^{*}H_1|_{\Delta_P})\cdots (f^{*}H_{d-1}|_{\Delta_P})\cdot (\mathcal{L}|_{\Delta_P})\right)&\\
         +\deg&(\Delta_P\to B)\left(H_1\cdots H_{d-1}\cdot A\right)\geqslant 0.
    \end{array}
\end{equation*}
Therefore,
\begin{equation*}
    \begin{array}{ll}
        h^{\mathcal{B}}_{(\mathcal{X},\mathcal{L})}(P) &=\dfrac{\deg\left((f^{*}H_1|_{\Delta_P})\cdots (f^{*}H_{d-1}|_{\Delta_P})\cdot (\mathcal{L}|_{\Delta_P})\right)}{\deg(\Delta_P\to B)}\\
          &\geqslant -\deg\left(H_1\cdots H_{d-1}\cdot A\right).
    \end{array}
\end{equation*}
\end{proof}

\begin{corollary}\label{uniform boundary of height}
If $\mathcal{L}_K\simeq\mathcal{O}_X$, then
 there exists a constant $C$ with
    $|h^{\mathcal{B}}_{(\mathcal{X},\mathcal{L})}(P)|\leq C $ for all $P\in X(\overline{K})$.
\end{corollary}

\begin{proof}
Apply Proposition \ref{upper boundary of height}
to $\mathcal{L}$ and $\mathcal{L}^{-1}$.
\end{proof}

\begin{proposition}
\label{uniform boundary of differnt height wrt same base}
Let $(\mathcal{X},\mathcal{L})$ and
$(\mathcal{X}',\mathcal{L}')$ be two $B$-models
of $(X,L)$.
Then there exist a constant $C>0$ with
\begin{equation*}
   | h^{\mathcal{B}}_{(\mathcal{X}',\mathcal{L}')}
    (P)-
    h^{\mathcal{B}}_{(\mathcal{X},\mathcal{L})}(P)|<C,
\end{equation*}
for all $P\in X(\overline{K})$.
\end{proposition}

\begin{proof}
Let $\mathcal{X}''$ be the graph of the birational map $\mathcal{X}'\to \mathcal{X}$ 
with projections
$\mu:\mathcal{X}''\to\mathcal{X}$ and
$\mu':\mathcal{X}''\to\mathcal{X}'$.
Let $
    \mathcal{M}
    :=
    (\mu')^*\mathcal{L}'
    \otimes(\mu^*\mathcal{L})^{-1}.$
We have $\mathcal{M}_K\simeq\mathcal{O}_X$ and $ h^{\mathcal{B}}_{(\mathcal{X}'',\mathcal{M})}(\cdot)= h^{\mathcal{B}}_{(\mathcal{X}'',(\mu')^*\mathcal{L}')}(\cdot)- h^{\mathcal{B}}_{(\mathcal{X}'',\mu^*\mathcal{L})}(\cdot)$, 
corollary \ref{uniform boundary of height}
completes the proof.
\end{proof}

\begin{remark}
We denote this height, defined up to a bounded function,
by $h_L^{\mathcal{B}}$ and call it the height associated
with $L$ and $\mathcal{B}$.
\end{remark}

\begin{proposition}\label{property of height}
Let $K/\mathbb{F}_q$ be a finitely generated extension
of transcendence degree $d\geq1$
and  $X$ be a  projective variety over $K$.
Fix a polarization
$\mathcal{B}=(B;H_1,\dots,H_{d-1})$ of $K$.
Then the following hold.
\begin{enumerate}
    \item[(1)]
    If $X\subseteq\mathbb{P}^n_K$ and
    $L=\mathcal{O}_{\mathbb{P}^n_K}(1)|_X$, then
    \begin{equation*}
        h_L^{\mathcal{B}}(\cdot)
        =
        h_{\mathrm{nv}}^{\mathcal{B}}|_{X(\overline{K})}(\cdot)
        +O(1).
    \end{equation*}

    \item[(2)]
    For any line bundles $L,L'$ on $X$, we have
    \begin{equation*}
        h_{L\otimes L'}^{\mathcal{B}}(\cdot)
        =
        h_L^{\mathcal{B}}(\cdot)+h_{L'}^{\mathcal{B}}(\cdot)+O(1).
    \end{equation*}
    In particular,
    \begin{equation*}
        h_{L^{-1}}^{\mathcal{B}}(\cdot)
        =
        -h_L^{\mathcal{B}}(\cdot)+O(1).
    \end{equation*}

    \item[(3)]
    The height $h_L^{\mathcal{B}}(\cdot)$ is bounded below on
    $(X\smallsetminus\SBs(L))(\overline{K})$, where
    \begin{equation*}
        \SBs(L):=\bigcap_{m\geq1}\Bs(L^{\otimes m}),
    \end{equation*}
    is the stable base locus of $L$.

    \item[(4)]
    Let $H_1',\dots,H_{d-1}'$ be nef line bundles on $B$
    such that $H_i'\otimes H_i^{-1}$ is
    $\mathbb{Q}$-effective for every $i$, and put
    \begin{equation*}
        \mathcal{B}'=(B;H_1',\dots,H_{d-1}').
    \end{equation*}
    If $L$ is ample, then
    \begin{equation*}
        h_L^{\mathcal{B}}(\cdot)
        \leq h_L^{\mathcal{B}'}(\cdot)+O(1).
    \end{equation*}
\end{enumerate}

\end{proposition}
\begin{proof}
  \item[(1)]By Proposition \ref{naive height given by moriwaki height},
the naive height is the model height associated with
$\mathcal{O}_{\mathbb{P}^n_B}(1)$.
Restricting to a model of $X$ and applying
Corollary \ref{uniform boundary of differnt height wrt same base},
we obtain the desired equality.

   \item[(2)]By definition, we have$h^{\mathcal{B}}_{(\mathcal{X},\mathcal{L}\otimes\mathcal{M})}=h^{\mathcal{B}}_{(\mathcal{X},\mathcal{L})}+h^{\mathcal{B}}_{(\mathcal{X},\mathcal{M})}$ and $h^{\mathcal{B}}_{(\mathcal{X},\mathcal{L}^{-1})}=-h^{\mathcal{B}}_{(\mathcal{X},\mathcal{L})}$. The assertions follow from the independence of model
   heights up to a bounded function.

 \item[(3)]By Noetherian, there is a positive integer $m$
such that
\begin{equation*}
    \SBs(L)=\Bs(L^{\otimes m}).
\end{equation*}
Applying Proposition \ref{upper boundary of height}
to $L^{\otimes m}$, we see that
$h_{L^{\otimes m}}^{\mathcal{B}}$ is bounded below on
$(X\smallsetminus\Bs(L^{\otimes m}))(\overline{K})$.
By (2),
\begin{equation*}
    h_{L^{\otimes(m)}}^{\mathcal{B}}
    =
    m h_L^{\mathcal{B}}+O(1).
\end{equation*}
Hence $h_L^{\mathcal{B}}$ is bounded below on
$(X\smallsetminus\SBs(L))(\overline{K})$.

\item[(4)]The assertion is immediate when $d=1$.
Assume that $d\geq2$.
By (2), replacing $L$ by a positive tensor power,
we may assume that $L$ is very ample.
Let $
    \phi_L:X\hookrightarrow\mathbb{P}^N_K$ 
be the corresponding embedding.
Let $\mathcal{X}$ be the  the closure
of $X$ in $\mathbb{P}^N_B$, and let
$f:\mathcal{X}\to B$ be the natural morphism.
Denote by $\mathcal{O}_{\mathcal{X}}(1)$ the pullback
of $\mathcal{O}_{\mathbb{P}^N_B}(1)$. Since $\mathcal{O}_{\mathcal{X}}(1)$ is $f$-ample,
there is an ample line bundle $\mathcal{M}$ on $B$
such that
\begin{equation*}
    \mathcal{L}
    :=
    \mathcal{O}_{\mathcal{X}}(1)\otimes f^*\mathcal{M}
\end{equation*}
is ample on $\mathcal{X}$.
Moreover, $\mathcal{L}_K\simeq L$.

For $P\in X(\overline{K})$, let $\Delta_P$ be the closure in $\mathcal{X}$.
Set $ A_i:=f^*H_i'|_{\Delta_P},
   $ , $
    B_i:=f^*H_i|_{\Delta_P},
    $ , $
    C:=\mathcal{L}|_{\Delta_P}.
$ We have
\begin{equation*}
        \begin{array}{rl}
            A_1\cdots A_{d-1}-&B_1\cdots B_{d-1}\\
            &= (A_1\cdots A_{d-1}- B_1\cdot A_2\cdots A_{d-1})+\cdots+(B_1\cdots B_{d-2}\cdot A_{d-1} -B_1\cdots B_{d-1})\\
            &=\sum\limits_{i=0}^{d-2}A_1\cdots A_i \cdot(A_{i+1}-B_{i+1})\cdot B_{i+2}\cdots B_{d-1}
        \end{array}
    \end{equation*}
Each $A_i-B_i$ is $\mathbb{Q}$-effective, while all
the other line bundles are nef. Therefore,
\begin{equation*}
    \deg\left(
        A_1\cdots A_{i-1}\cdot(A_i-B_i)
        \cdot B_{i+1}\cdots B_{d-1}\cdot C
    \right)\geq0.
\end{equation*}
Summing over $i$, we obtain
\begin{equation*}
    \deg\left(A_1\cdots A_{d-1}\cdot C\right)
    \geq
    \left(B_1\cdots B_{d-1}\cdot C\right).
\end{equation*}
Which gives that
\begin{equation*}
    h^{\mathcal{B}}_{(\mathcal{X},\mathcal{L})}(P)
    \leq
    h^{\mathcal{B}'}_{(\mathcal{X},\mathcal{L})}(P).
\end{equation*}
Hence, we obatin (4).
\end{proof}

We now discuss the Northcott property. We first consider
a special case.

\begin{lemma}\label{self interestion number of line bundle}
Let $B=(\mathbb{P}^1_{\mathbb{F}_q})^d$, where $d\geq1$,
and let $p_i$ be the $i$-th projection.
Let $L_1,\dots,L_d$ be line bundles on
$\mathbb{P}^1_{\mathbb{F}_q}$, and set
\begin{equation*}
    L:=p_1^*L_1\otimes\cdots\otimes p_d^*L_d,
    \qquad
    \Delta_\infty^i:=p_i^*(\infty).
\end{equation*}
Then
\begin{equation*}
    \deg\left((L|_{\Delta_\infty^i})^{d-1}\right)
    =(d-1)!
    \prod_{\substack{1\leq j\leq d\\j\neq i}}
    \deg(L_j).
\end{equation*}
\end{lemma}
\begin{proof}
The case $d=1$ is immediate.
Assume that $d\geq2$.
After relabeling, we may assume that $i=d$.
There is an isomorphism
 \begin{equation}
        \phi:\underbrace{\mathbb{P}^1_{{\mathbb{F}_q}}\times_{\mathbb{F}_q}\cdots \times_{\mathbb{F}_q}\mathbb{P}^1_{{\mathbb{F}_q}}
        }_{d-1} \to \Delta^d_{\infty} 
     \end{equation}
Let $q_j$ be the
$j$-th projection.
For $1\leq j\leq d-1$, we have
\begin{equation*}
    \phi^*(p_j^*L_j|_{\Delta_\infty^d})\simeq q_j^*L_j,
\end{equation*}
while
$\phi^*(p_d^*L_d)\simeq\mathcal{O}_Y$.
Thus
\begin{equation*}
    \phi^*(L|_{\Delta_\infty^d})
    \simeq q_1^*L_1\otimes\cdots\otimes q_{d-1}^*L_{d-1}.
\end{equation*}
Now we using it to calculate the interestion number
\begin{equation*}
\begin{aligned}
   \deg \left((L|_{\Delta_\infty^d})^{d-1}\right)
    &=
    \deg\left(
        (q_1^*L_1+\cdots+q_{d-1}^*L_{d-1})^{d-1}
    \right)\\
    &=
    \sum_{a_1+\cdots+a_{d-1}=d-1}
    \frac{(d-1)!}{a_1!\cdots a_{d-1}!}
    \deg\left(
        (q_1^*L_1)^{a_1}\cdots
        (q_{d-1}^*L_{d-1})^{a_{d-1}}
    \right).
\end{aligned}
\end{equation*}
Since each $q_j^*L_j$ is pulled back from $\mathbb{P}^1_{\mathbb{F}_q}$,
every term with $a_j\geq2$ vanishes.
Hence only the term $a_1=\cdots=a_{d-1}=1$ can contribute.
Therefore,
\begin{equation*}
    \deg\left((L|_{\Delta_\infty^d})^{d-1}\right)
    =(d-1)!
   \deg \left(q_1^*L_1\cdots q_{d-1}^*L_{d-1}\right)
    =(d-1)!\prod_{j=1}^{d-1}\deg(L_j).
\end{equation*}
\end{proof}
\begin{lemma}\label{Northcott of special case }
Let
\begin{equation*}
    K=\mathbb{F}_q(z_1,\dots,z_d),
    \qquad
    B=(\mathbb{P}^1_{\mathbb{F}_q})^d,
    \qquad
    H=\bigotimes_{i=1}^d p_i^*\mathcal{O}_{\mathbb{P}^1}(1),
\end{equation*}
and put $\mathcal{B}=(B;H,\dots,H)$ with $d-1$ copies
of $H$.
Then, for every real number $M$, the set
\begin{equation*}
    \left\{
        P\in\mathbb{P}^n(K)
        \mid h_{\mathrm{nv}}^{\mathcal{B}}(P)\leq M
    \right\}
\end{equation*}
is finite.
\end{lemma}

\begin{proof}
For $1\leq i\leq d$, set
\begin{equation*}
    A_i:=
    \bigotimes_{\substack{1\leq j\leq d\\j\neq i}}
    p_j^*\mathcal{O}_{\mathbb{P}^1}(1),
    \qquad
    \mathcal{B}_i:=(B;A_i,\dots,A_i).
\end{equation*}
By Lemma \ref{self interestion number of line bundle}, we have
\begin{equation*}
    \left(A_i^{d-1}\cdot\Delta_\infty^j\right)
    =
    \begin{cases}
        (d-1)!,& i=j,\\
        0,& i\neq j.
    \end{cases}
\end{equation*}
Since $H\otimes A_i^{-1}=p_i^*\mathcal{O}_{\mathbb{P}^1}(1)$
is effective, Proposition \ref{property of height}(4)
gives a constant $C_i$ such that
\begin{equation*}
    h_{\mathrm{nv}}^{\mathcal{B}_i}(P)
    \leq h_{\mathrm{nv}}^{\mathcal{B}}(P)+C_i.
\end{equation*}

Write $P=[f_0:\cdots:f_n]$, where
$f_0,\dots,f_n\in\mathbb{F}_q[z_1,\dots,z_d]$
with $(f_0,\dots,f_n)=1$.
For every prime divisor $\Gamma$ meeting
$\mathbb{A}^d_{\mathbb{F}_q}$, we have
\begin{equation*}
    \max_{0\leq\ell\leq n}
    \{-\ord_\Gamma(f_\ell)\}=0.
\end{equation*}
At infinity,
\begin{equation*}
    \ord_{\Delta_\infty^j}(f_\ell)
    =-\deg_{z_j}(f_\ell).
\end{equation*}
Therefore,
\begin{equation*}
\begin{aligned}
    h_{\mathrm{nv}}^{\mathcal{B}_i}(P)
    &=
    \sum_\Gamma
    \max_{0\leq\ell\leq n}\{-\ord_\Gamma(f_\ell)\}
    \left(\Gamma\cdot A_i^{d-1}\right)\\
    &=
    \sum_{j=1}^d
    \max_{0\leq\ell\leq n}\{\deg_{z_j}(f_\ell)\}
    \left(\Delta_\infty^j\cdot A_i^{d-1}\right)\\
    &=
    (d-1)!\max_{0\leq\ell\leq n}
    \{\deg_{z_i}(f_\ell)\}.
\end{aligned}
\end{equation*}
If $h_{\mathrm{nv}}^{\mathcal{B}}(P)\leq M$, then
\begin{equation*}
    (d-1)!\max_{0\leq\ell\leq n}\{\deg_{z_i}(f_\ell)\}\leq M+C_i
\end{equation*}
for every $i$.
Thus the degree of each $f_i$ in each variable is
bounded independently of $P$.
Since $\mathbb{F}_q$ is finite, there are only finitely
many such polynomials, and hence only finitely many
such points.
\end{proof}

We now prove the Northcott property for a general
finitely generated field.
\begin{theorem}\label{Northcott}
Let $K/\mathbb{F}_q$ be a finitely generated extension
of transcendence degree $d\geq1$, and let
$\mathcal{B}=(B;H_1,\dots,H_{d-1})$ be a big polarization.
Let $X$ be a projective variety over $K$, and let $L$
be an ample line bundle on $X$.
Then, for every real number $M$ and every positive
integer $e$, the set
\begin{equation*}
    \left\{
        P\in X(\overline{K})
        \mid h_L^{\mathcal{B}}(P)\leq M,\quad
        [K(P):K]\leq e
    \right\}
\end{equation*}
is finite.
\end{theorem}
Set
\begin{equation*}
    K_0=\mathbb{F}_q(z_1,\dots,z_d),
    \qquad
    B_0=(\mathbb{P}^1_{\mathbb{F}_q})^d,
    \qquad
    H_0=\bigotimes_{i=1}^d p_i^*\mathcal{O}_{\mathbb{P}^1}(1),
\end{equation*}
and put $\mathcal{B}_0=(B_0;H_0,\dots,H_0)$.
We prove the theorem through the following reductions.

\begin{lemma}\label{special case 1}
Theorem \ref{Northcott} holds for $\mathcal{B}_0$
and $e=1$.
\end{lemma}

\begin{proof}
Since $L$ is ample, we can take an integer $m$ larger enough, such that
$L^{\otimes m}$ defines an embedding
$
    \phi_{L^{\otimes m}}:X\hookrightarrow\mathbb{P}^N_{K_0}.
$
By Proposition \ref{property of height}(1) and (2),
\begin{equation*}
    h_{\mathrm{nv}}^{\mathcal{B}_0}
    \left(\phi_{L^{\otimes m}}(P)\right)
    =
    m h_L^{\mathcal{B}_0}(P)+O(1).
\end{equation*}
Thus a bound on $h_L^{\mathcal{B}_0}(P)$ gives a bound
on the naive height.
The result follows from
Lemma \ref{Northcott of special case }.
\end{proof}
We next reduce the general polarization to $\mathcal{B}_0$.

\begin{lemma}\label{special case 2}
If Theorem \ref{Northcott} holds for $\mathcal{B}_0$,
then it holds for every big polarization $\mathcal{B}$.
\end{lemma}

\begin{proof}
    As in the previous lemma, it suffices to consider
$(X,L)=(\mathbb{P}^N_K,\mathcal{O}(1))$ and use the naive
height.
Choose a transcendence basis $z_1,\dots,z_d$ of
$K/\mathbb{F}_q$, set
$K_0=\mathbb{F}_q(z_1,\dots,z_d)$ and $m=[K:K_0]$.
The inclusion $K_0\subseteq K$ induces a rational map
$B\dashrightarrow B_0$. Replacing $B$ by the normalization of the closure of its graph,
 we may assume there
is a dominant generically finite morphism $\mu:B\to B_0$. 
Since each $H_i$ is big, there is a positive integer
$n_i$ such that
$H_i^{\otimes n_i}\otimes(\mu^*H_0)^{-1}$ is effective.
By Proposition \ref{height under finite extension},
and Proposition \ref{property of height}(4), we obtain
\begin{equation*}
    m h_{\mathrm{nv},K_0}^{\mathcal{B}_0}(P)=h_{L}^{(B;\mu^*H_0,\dots,\mu^*H_0)}(P)\leq a h_{L}^{\mathcal{B}}(P)+b
\end{equation*}
for some constants $a,b>0$. The Northcott property for $\mathcal{B}_0$ therefore
gives the required finiteness.
\end{proof}

\begin{lemma}\label{special case 3}
If Theorem \ref{Northcott} holds for $e=1$,
then it holds for every positive integer $e$.
\end{lemma}

\begin{proof}
For each positive integer $e$, we prove
that the set
\begin{equation*}
    \left\{
        P\in\mathbb{P}^N(\overline{K})
        \mid h_{\mathrm{nv},K}^{\mathcal{B}}(P)\leq M,
        \ [K(P):K]=e
    \right\}
\end{equation*}
is finite.

Let $P$ belong to this set and   $\{\sigma_1,\dots,\sigma_e\}$ be the set of all embeddings of $K(P)$ into $\overline{K}$ over $K$ and 
write $P_i$ the point given by $\Spec(\overline{K})\xrightarrow{\sigma_i^{*}}\Spec(K(P))\xrightarrow{P}X$ counted with multiplicity. 
Let $\mathcal{Y}$ be the main component of $\underbrace{\mathcal{X}\times_B \cdots\times_{B}\mathcal{X}}_{e} $
 and $\mathcal{M}=p_1^*\mathcal{L}\otimes\cdots\otimes p_e^*\mathcal{L}$,
where $p_i$ is the $i$-th projection.
Denote the closures of $P_i$ in $\mathcal{X}$ by $\Delta_{P_i}$
Since $\Delta_{P_i}=\Delta_P$ for every $i$,
the projection formula gives
\begin{equation*}
    \begin{array}{rl}
        h_{\left(\mathcal{Y},\mathcal{M}\right)}^{\mathcal{B}}&(P_1,\dots,Pe)\\
       & =\dfrac{\left(\left(f^{\prime}\right)^{*}\left(H_1\cdots H_{d-1}\right)\cdot\mathcal{M}\cdot \Delta_{(P_1,\dots,P_e)}\right)}{\deg\left(\Delta_{P_1,\dots,P_e}\to B\right)}\\
       &=\dfrac{\sum\limits_{i=1}^{e}\left(p_i^{*}f^{*}\left(H_1\cdots H_{d-1}\right)\cdot p_i^{*}(\mathcal{L})\cdot \Delta_{(P_1,\dots,P_e)} \right)}{\deg\left(\Delta_{P_1,\dots,P_e}\to B\right)}\\
       &=\dfrac{\sum\limits_{i=1}^{e}\deg\left(\Delta_{(P_1,\dots,P_e)}\to \Delta_{P_i}\right)\left(f^{*}\left(H_1\cdots H_{d-1}\right)\cdot \mathcal{L}\cdot \Delta_{P_i} \right)}{\deg\left(\Delta_{P_1,\dots,P_e}\to B\right)}\\
       &=eh^{\mathcal{B}}_{(\mathcal{X},\mathcal{L})}(P)
    \end{array}
\end{equation*}
where $f:\mathcal{X}\to B$ and $f':\mathcal{Y}\to B$ are canonical morphism.

Let $\rho:X^e\to Z=X^e/\mathfrak{S}_e$ be the quotient
morphism, where $\mathfrak{S}_e$ is the symmetric group.
Since 
$M=\mathcal{M}_K=p_1^*L\otimes\cdots\otimes p_e^*L$ is invariant under the action of $\mathfrak{S}_e$ it will descends naturally to a line bundle $N$ on
$Z=X^e/\mathfrak{S}_e$, with $\rho^*N\simeq M$.
By functoriality of heights, we have
\begin{equation*}
    h_N^{\mathcal{B}}\circ\rho
    =
    h_{(\mathcal{Y},\mathcal{M})}^{\mathcal{B}}+O(1).
\end{equation*}
Hence, there exist a constant $C$ independent of $P$, such that
\begin{equation*}
    h_N^{\mathcal{B}}(\rho(P_1,\dots,P_e))\leq eh_L^{\mathcal{B}}(P)+C\leq eM+C
\end{equation*}
Moreover, Put $Y=\Spec K(P)$
,the closed subscheme $Y\subset X$ is finite flat
of degree $e$ over $K$, its canonical section $t\in\operatorname{Sym}^e_K(Y)(K)$ induces,
via $Y\hookrightarrow X$, a point $z\in\operatorname{Sym}^e_K(X)(K)$
\cite[ (6.3.4.2)]{SGA4}. we then get  $\rho(P_1,\dots,P_e)\in Z(K)$.
The Northcott property
for $e=1$ gives that there are only finitely many
$\rho(P_1,\dots,P_e)$ with
\begin{equation*}
    h_N^{\mathcal{B}}(\rho(P_1,\dots,P_e))\leq eM+C
\end{equation*}
But fiber of $\rho$ contains at most $e!$ points.
there are only finitely many
$(P_1,\dots,P_e)$, and therefore only finitely many
such points $P$.
This completes the proof.

\end{proof}

\begin{proof}[Proof of Theorem \ref{Northcott}]
By Lemma \ref{special case 1}, the Northcott property
holds for $\mathcal{B}_0$ when $e=1$.
Lemma \ref{special case 3} extends it to every positive
integer $e$, and Lemma \ref{special case 2} then gives
the result for every big polarization $\mathcal{B}$.
\end{proof}

\subsection{Comparison of heights under field extensions}
In this section, we study the comparability of Moriwaki
heights under finitely generated field extensions.
Let $k=\mathbb{F}_q$ and $K/k$ be a finitely generated
extension with $\trdeg_k(K)=d\geq1$. Let $X$ be a
projective variety over $K$.
Fix an ample line bundle $L$ on $X$, a big polarization
$\mathcal{B}$ of $K$, and a model
$(\mathcal{X},\mathcal{L})$ of $(X,L)$.

Let $K'/K$ be an extension of finitely generated fields
over $k$.
Fix a $K$-embedding $\overline{K}\hookrightarrow\overline{K'}$.
For $x\in X(\overline{K})$, denote by
$x_{K'}\in X_{K'}(\overline{K'})$ the point obtained
by base change.
Choose a big polarization $\mathcal{B}'$ of $K'$ and
a model $(\mathcal{X}',\mathcal{L}')$ of
$(X_{K'},L')$, where $L'$ is ample.
The corresponding Moriwaki height then restricts to a
function on $X(\overline{K})$.
We will compare this function with
$h^{\mathcal{B}}_{(\mathcal{X},\mathcal{L})}$.

\begin{definition}
With the notation above, we say that
$h^{\mathcal{B}}_{(\mathcal{X},\mathcal{L})}$ and
$h^{\mathcal{B}'}_{(\mathcal{X}',\mathcal{L}')}$ are
\textbf{comparable} if there exist constants
$N_1,N_2>0$ and $M_1,M_2\in\mathbb{R}$ such that
\begin{equation*}
    N_1h^{\mathcal{B}'}_{(\mathcal{X}',\mathcal{L}')}(x_{K'})
    +M_1
    \leq h^{\mathcal{B}}_{(\mathcal{X},\mathcal{L})}(x)\leq
    N_2h^{\mathcal{B}'}_{(\mathcal{X}',\mathcal{L}')}(x_{K'})
    +M_2
\end{equation*}
for every $x\in X(\overline{K})$.
\end{definition}

The following transitivity property follows directly
from the definition; see also \cite{LS2025}.
\begin{lemma}
Let $K_1/K$ and $K_2/K$ be extensions of finitely
generated fields over $k$.
Let
$h(\cdot)=h^{\mathcal{B}}_{(\mathcal{X},\mathcal{L})}(\cdot)$ and
$h_i(\cdot)=h^{\mathcal{B}_i}_{(\mathcal{X}_i,\mathcal{L}_i)}(\cdot)$
for $i=1,2$ be the corresponding Moriwaki heights,
viewed as functions on $X(\overline{K})$.
Consider the following conditions:
\begin{itemize}
    \item[(a)] $h(\cdot)$ and $h_1(\cdot)$ are comparable.
    \item[(b)] $h(\cdot)$ and $h_2(\cdot)$ are comparable.
    \item[(c)] $h_1(\cdot)$ and $h_2(\cdot)$ are comparable.
\end{itemize}
If (a) and (c) hold, then (b) holds.
Similarly, if (b) and (c) hold, then (a) holds.
We first compare heights over a fixed field.
\end{lemma} 
\begin{lemma}\label{comparable in same field}
Let $\mathcal{B}=(B;H_1,\dots,H_{d-1})$ and
$\mathcal{B}'=(B';H_1',\dots,H_{d-1}')$ be two big
polarizations of $K$.
Let $X$ be a projective variety over $K$ and 
$L,L'$ be ample line bundles on $X$.
Choose models $(\mathcal{X},\mathcal{L})$ and
$(\mathcal{X}',\mathcal{L}')$ of $(X,L)$ and
$(X,L')$, respectively.
Then
$h^{\mathcal{B}}_{(\mathcal{X},\mathcal{L})}(\cdot)$ and
$h^{\mathcal{B}'}_{(\mathcal{X}',\mathcal{L}')}(\cdot)$ are
comparable.
\end{lemma}

\begin{proof}
Let $B''$ be the normalization of the graph of the
birational map $B'\dashrightarrow B$, and let
$g:B''\to B$ and $g':B''\to B'$ be the projections.
Set
\begin{equation*}
    \widetilde{\mathcal{B}}
    :=(B'';g^*H_1,\dots,g^*H_{d-1}),\quad
    \widetilde{\mathcal{B}}'
    :=(B'';g'^*H_1',\dots,g'^*H_{d-1}').
\end{equation*}
Since $g$ and $g'$ are birational, all the line bundles
$g^*H_i$ and $g'^*H_i'$ are nef and big.
By Kodaira's Lemma \cite[Proposition~10.10]{Vojta2009}, we can choose a positive integer $N$ such that
both
$(g^*H_i)^{\otimes N}\otimes(g'^*H_i')^{-1}$ and
$(g'^*H_i')^{\otimes N}\otimes(g^*H_i)^{-1}$
are effective for every $i$.
By Proposition \ref{property of height}(4), we have
\begin{equation*}
        h^{(B^{''},g^{\prime*}H_1^{\prime},\dots,g^{\prime*}H_{d-1}^{\prime})}_{(\mathcal{X}_{B^{''}},\mathcal{L}_{B^{''}})}(\cdot)\leq N^{d-1} h^{(B^{''},g^{*}H_1,\dots,g^{*}H_{d-1})}_{(\mathcal{X}_{B^{''}},\mathcal{L}_{B^{''}})}(\cdot)+O(1)
\end{equation*}
and similarly,
\begin{equation*}
     h^{(B^{''},g^{*}H_1,\dots,g^{*}H_{d-1})}_{(\mathcal{X}_{B^{''}},\mathcal{L}_{B^{''}})}(\cdot)  \leq N^{d-1} h^{(B^{''},g^{\prime*}H_1^{\prime},\dots,g^{\prime*}H_{d-1}^{\prime})}_{(\mathcal{X}_{B^{''}},\mathcal{L}_{B^{''}})}(\cdot)+O(1)
\end{equation*}
Hence these two heights are comparable.
On the other hand, since $g$ and $g'$ are birational,
the projection formula gives
\begin{equation*}
    h^{\widetilde{\mathcal{B}}}
    _{(\mathcal{X}_{B''},\mathcal{L}_{B''})}(\cdot)
    =
    h^{\mathcal{B}}_{(\mathcal{X},\mathcal{L})}(\cdot),\quad
    h^{\widetilde{\mathcal{B}}'}
    _{(\mathcal{X}'_{B''},\mathcal{L}'_{B''})}(\cdot)
    =
    h^{\mathcal{B}'}_{(\mathcal{X}',\mathcal{L}')}(\cdot).
\end{equation*}

Since $L$ and $L'$ are ample, we can choose positive
integers $a,b$ such that
$L^{\otimes a}\otimes(L')^{-1}$ and
$(L')^{\otimes b}\otimes L^{-1}$ are ample.
Their heights with respect to
$\widetilde{\mathcal{B}}'$ are bounded below.
So we obtain
\begin{equation*}
    h^{\widetilde{\mathcal{B}}'}
    _{(\mathcal{X}'_{B''},\mathcal{L}'_{B''})}(\cdot)
    \leq
    a\,h^{\widetilde{\mathcal{B}}'}
    _{(\mathcal{X}_{B''},\mathcal{L}_{B''})}(\cdot)+O(1),\quad
    h^{\widetilde{\mathcal{B}}'}
    _{(\mathcal{X}_{B''},\mathcal{L}_{B''})}(\cdot)
    \leq
    b\,h^{\widetilde{\mathcal{B}}'}
    _{(\mathcal{X}'_{B''},\mathcal{L}'_{B''})}(\cdot)+O(1).
\end{equation*}
Thus these two heights are comparable.
Combining this with the preceding comparison and
the projection formula proves the lemma.
\end{proof}

We now consider finitely generated field extensions.

\begin{theorem}\label{comparable of  height under field extension}
Let $K_1/K$ be an extension of finitely generated
fields over $k$, $\mathcal{B}$ and $\mathcal{B}_1$ be big
polarizations of $K$ and $K_1$, respectively.
Choose models $(\mathcal{X},\mathcal{L})$ and
$(\mathcal{X}_1,\mathcal{L}_1)$ of $(X,L)$ and
$(X_{K_1},L_1)$, where $L$ and $L_1$ are ample.
Then there exist constants $N_1,N_2>0$ and
$M_1,M_2\in\mathbb{R}$ such that
\begin{equation*}
    N_1h^{\mathcal{B}_1}_{(\mathcal{X}_1,\mathcal{L}_1)}
    (x_{K_1})+M_1
    \leq h^{\mathcal{B}}_{(\mathcal{X},\mathcal{L})}(x)
    \leq
    N_2h^{\mathcal{B}_1}_{(\mathcal{X}_1,\mathcal{L}_1)}
    (x_{K_1})+M_2
\end{equation*}
for every $x\in X(\overline{K})$.
\end{theorem}

\begin{proof}
Choose a transcendence basis $t_1,\dots,t_e$ of
$K_1/K$, so that $K_1/K(t_1,\dots,t_e)$ is finite.
By Proposition \ref{height under finite extension},
the height under a finite extension  up to
the degree of the extension.
Together with Lemma \ref{comparable in same field},
this reduces the proof to $K(t_1,\dots,t_e)/K$.
And it suffices to consider $K(t)/K$.
Moreover, Lemma \ref{comparable in same field}
allows us to choose convenient polarizations and
models.

Write $d=\trdeg_k(K)$.
Choose an ample line bundle $H$ on $B$ and put
$\mathcal{B}=(B;H,\dots,H)$, with $d-1$ copies of $H$.
Let $B'=\mathbb{P}^1_k\times_k B$ and let $p_1,p_2$
be the projections.
Set
$H'=p_1^*\mathcal{O}(1)\otimes p_2^*H$ and
$\mathcal{B}'=(B';H',\dots,H')$.
Denote $\mathcal{X}'=\mathcal{X}\times_B B'$,  $\pi:\mathcal{X}'\to\mathcal{X}$ be the projection
and $\mathcal{L}'=\pi^*\mathcal{L}$.
Consider the following diagram
\begin{center}
\begin{tikzcd}
    \mathcal{X}'
        \arrow[r,"f'"]
        \arrow[d,"\pi"']
    & B'
        \arrow[r,"p_1"]
        \arrow[d,"p_2"]
    & \mathbb{P}^1_k
        \arrow[d]\\
    \mathcal{X}
        \arrow[r,"f"]
    & B
        \arrow[r]
    & \Spec k .
\end{tikzcd}
\end{center}
Write $q_1=p_1\circ f'$.

For $P\in X(\overline{K})$, let $\Delta_P$ be its
closure in $\mathcal{X}$, and let $\Delta_P'$ be
the closure of $P_{K(t)}$ in $\mathcal{X}'$.
We have
\begin{equation*}
    \begin{array}{rl}
           h^{\mathcal{B}^{\prime}}_{(\mathcal{X}^{\prime},\mathcal{L}^{\prime})}(P_{K(t)})&=\dfrac{\left(\left(\pi^{*}\mathcal{L}|_{\Delta_P^{\prime}}\right)\cdot  \left(f^{\prime*}\left(p_1^{*}\mathcal{O}(1)\otimes p_2^{*}H\right)|_{\Delta_P^{'}}\right)^{d}\right)}{[K(P):K]}\\
    \end{array}
\end{equation*}
Using additive notation in intersection products,
we obtain
\begin{equation*}
\begin{aligned}
    \left(f^{\prime*}\left(p_1^{*}\mathcal{O}(1)\otimes p_2^{*}H\right)\right)^{d}
    &=
    \sum_{j=0}^{d}
    \binom{d}{j}
    \left(q_1^*\mathcal{O}(1)\right)^j
    \cdot(\pi^*f^*H)^{d-j}\\
    &=
    (\pi^*f^*H)^d
    +d\,q_1^*\mathcal{O}(1)
    \cdot(\pi^*f^*H)^{d-1},
\end{aligned}
\end{equation*}
since every term with $j\geq2$ vanishes.

Under the natural identification
$\mathcal{X}'\simeq\mathcal{X}\times_k\mathbb{P}^1_k$,
we have $\Delta_P'=\Delta_P\times_k\mathbb{P}^1_k$.
Since $\dim\Delta_P=d$, it follows that
\begin{equation*}
    \left(
        \pi^*\mathcal{L}\cdot(\pi^*f^*H)^d
        \cdot\Delta_P'
    \right)=0.
\end{equation*}
For the second term, we have
\begin{equation*}
\begin{aligned}
    \left(
        \pi^*\mathcal{L}
        \cdot q_1^*\mathcal{O}(1)
        \cdot(\pi^*f^*H)^{d-1}
        \cdot\Delta_P'
    \right)
    &=
    \deg_{\mathbb{P}^1_k}\mathcal{O}(1)
    \left(
        \mathcal{L}\cdot(f^*H)^{d-1}\cdot\Delta_P
    \right)\\
    &=
    \left(
        \mathcal{L}\cdot(f^*H)^{d-1}\cdot\Delta_P
    \right).
\end{aligned}
\end{equation*}
Consequently,
\begin{equation*}
    h^{\mathcal{B}'}_{(\mathcal{X}',\mathcal{L}')}
    (P_{K(t)})
    =
    d\,h^{\mathcal{B}}_{(\mathcal{X},\mathcal{L})}(P),
\end{equation*}
for every $P\in X(\overline{K})$.
Combining this identity with
Lemma \ref{comparable in same field} and
Proposition \ref{height under finite extension}
completes the proof.
\end{proof}

\begin{proposition}\label{height rela }
Let $X$ be a normal projective variety over $K$.  Let $D$ be a numerically trivial $\mathbb{R}$-Cartier
divisor and $A$ an ample Cartier divisor on $X$.
Fix associated heights $h_D$ and $h_A$ with respect
to a fixed polarization, and assume that $h_A\geq1$.
Then there exists a constant $C>0$ such that
\begin{equation*}
    h_D(x)\leq C\sqrt{h_A(x)}
\end{equation*}
for every $x\in X(\overline{K})$.
\end{proposition}
\begin{proof}
By $\mathbb{R}$-linearity of heights, we may assume $D$ and $A$ are Cartier divisors. 
 Since $X$ is a normal projective variety,
after a finite purely inseparable extension and
normalization, followed by a finite separable extension,
we obtain a finite extension $F/K$ and a finite
surjective morphism
\begin{equation*}
    \nu:\coprod_{j=1}^{s}Y_j\longrightarrow X_F,
\end{equation*}
where each $Y_j$ is geometrically normal,
geometrically integral and projective; see
\cite[Tag 0BXT, Tag 054R]{stacks-project}.
Write $\nu_j=\nu|_{Y_j}$, and set
$D_j=\nu_j^*D_F$ and $A_j=\nu_j^*A_F$, we have $D_j$ is numerically trivial and $A_j$ is ample. Replacing $X$ by $Y_j$, by functoriality of heights. We may assume $X$ is     geometrically normal and geometrically integral. Furthermore, since $G_n/G_a$
is finite by \cite{Mats1957}, there exists
a positive integer $m$ such that $mD$ is algebraically
trivial. So we may assume that $D$ is algebraically trivial.

Since those property are stable under base change, we may further
replace $K$ by a finite extension and assume that
$X$ has a $K$-rational point.
Fixing such a point, let $u_X:X\to\underline{\Pic}^0(P(X))=\Alb(X)$
be the Albanese morphism with $\Alb(X)$ projective, see \cite[Theorem~9.5.4]{FAG2005}.
There exists an
algebraically trivial divisor $E$ on $\Alb(X)$
such that $D\sim u_X^*E$. Choose a symmetric ample line bundle $E_1$ on $\Alb(X)$.
Then the morphism
\begin{equation*}
    \phi:\Alb(X)\longrightarrow\Pic^0(\Alb(X)),
    \qquad
    a\longmapsto t_a^*E_1-E_1
\end{equation*}
is surjective. Hence there exists
$c\in\Alb(X)(\overline{K})$ such that
$E\sim t_c^*E_1-E_1$.
    For $P\in X(\overline{K})$, put $Q=u_X(P)$.
By functoriality and additivity of heights, we have
\begin{equation*}
    \begin{array}{rl}
        h_{X,D}(P)&=h_{X,u_X^{*}E}(P)+O(1)\\
        &=h_{\Alb(X),E}(Q)+O(1)\\
        &=\widehat{h}_{\Alb(X),E}(Q)+O(1)\\
        &=\widehat{h}_{\Alb(X),t_c^{*}E_1}(Q)-\widehat{h}_{\Alb(X),E_1}(Q)+O(1)\\
        &=\widehat{h}_{\Alb(X),E_1}(t_cQ)-\widehat{h}_{\Alb(X),E_1}(t_cO)-\widehat{h}_{\Alb(X),E_1}(Q)+O(1)\\
        &=\widehat{h}_{\Alb(X),E_1}(Q+c)-\widehat{h}_{\Alb(X),E_1}(c)-\widehat{h}_{\Alb(X),E_1}(Q)+O(1)\\
        &=2\langle Q,c \rangle_{E_1}+O(1)\\
        &\leq2\sqrt{\widehat{h}_{\Alb(X),E_1}(Q)\widehat{h}_{\Alb(X),E_1}(c)}+O(1).
    \end{array}
\end{equation*}
Similarly, we have
\begin{equation*}
    \begin{array}{rl}
        \widehat{h}_{\Alb(X),E_1}(Q)&=   h_{\Alb(X),E_1}(u_X( P ))+O(1)\\
        &=h_{X,u_X^{*}E_1}(P)+O(1)\\
        &\leq M h_{X,A}(P)+O(1).
        \end{array}
\end{equation*}
On the other hand, the $O(1)$ terms are uniformly
bounded above by a constant $N$.
Since $h_A(x)\geq1$, we may choose a sufficiently
large constant $C>0$ such that
\begin{equation*}
    h_D(x)\leq C\sqrt{h_A(x)},
\end{equation*}
for every $x\in X(\overline{K})$.
\end{proof}

\section{Dynamics over finitely generated fields of positive characteristic}\label{fg case chatper}

\subsection{Some intersection inequalities in positive characteristic}
In this section, we fix a finite field $k$.  
 Let $K/k$ be a finitely generated field with $\trdeg_k(K)=d\geq 1$, take a polarization $B$ of $K$. Let $X$ be a normal projective variety over $K$, we take a $B$-model $\mathcal{X}$ . For a cycle class $\mathcal{Z}\in \ch^{j}(\mathcal{X})$, we say $\mathcal{Z}\geq 0$ if for any nef line bundles $\mathcal{L}_1$, $\dots$, $\mathcal{L}_{n-j}$, the interestion number
  \begin{equation*}
    \left(\mathcal{L}_1\cdots \mathcal{L}_{n-k}\cdot  \mathcal{Z}\right)\geq 0.
  \end{equation*} 
  We first establish an intersection inequality for nef
line bundles, which will be used in the height estimates
below. The case $r=0$ follows from the bigness 
inequality, and the general case is proved by induction
on $r$.
  \begin{lemma}\label{num inequality}
    Let $\mathcal{X}\to k$ be a normal projective variety with dimension $n$, and $0\leq r\leq n-1$  be an integer. Let $\mathcal{L}$, $\mathcal{M}$, $\mathcal{N}_i$, $i=1,2,\dots, r$ be nef line bundles such that 
      \begin{equation*}
      \left(\mathcal{M}^{n-r}\cdot \mathcal{N}_1\cdots \mathcal{N}_r\right)>\left(n-r\right) \left(\mathcal{M}^{n-r-1}\cdot \mathcal{L}\cdot \mathcal{N}_1\cdots \mathcal{N}_r\right).
  \end{equation*}
Then we have 
    \begin{equation*}
      \mathcal{L}\cdot \mathcal{N}_1\cdots \mathcal{N}_r\leq \mathcal{M}\cdot \mathcal{N}_1\cdots \mathcal{N}_r.
  \end{equation*}
  \end{lemma}
 
\begin{proof}
    Sine the interestion numbers are invariant under field extension, see \cite[Example~6.2.9]{Ful98}. We may replace $k$ by $k(t)$, and assume the field $k$ is infinite field. Thus we can use the Bertini theorem for infinite field field. We prove this inequality by induction on r   .
When $r=0$, 
\cite[Proposition 5.4]{Cutk15} gives that
\begin{equation*}
    \vol
    (\mathcal{M}-\mathcal{L})
    \geq
    \left(\mathcal{M}^{n}\right)
    -
    n\left(\mathcal{M}^{n-1}\cdot\mathcal{L}\right)
    >0.
\end{equation*}
Thus $\mathcal{M}-\mathcal{L}$ is big. Consequently, there
exist a positive integer $a$ and an effective line bundle $\mathcal{E}$ such that
\begin{equation*}
   \mathcal{E}\equiv a(\mathcal{M}-\mathcal{L}).
\end{equation*}
For any nef line bundles
$\mathcal{H}_1,\dots,\mathcal{H}_{n-1}$, we obtain
\begin{equation*}
    \left(
    (\mathcal{M}-\mathcal{L})\cdot
    \mathcal{H}_1\cdots\mathcal{H}_{n-1}
    \right)
    =
    \frac{1}{a}
    \left(
    \mathcal{E}\cdot\mathcal{H}_1\cdots\mathcal{H}_{n-1}
    \right)
    \geq 0.
\end{equation*}
This proves the case $r=0$.

Assume that $r\geq1$ and that the statement holds for
$r-1$.
since nef line bundles can be approximated by ample line bundles, we can
assume all these line bundles are ample. 
Choose a positive integer $m$ such that
$\mathcal{N}_1^{\otimes m}$ is very ample, and set
$\Gamma=|\mathcal{N}_1^{\otimes m}|$. Since $k$ is infinite, by Seidenberg's Bertini theorem
\cite[Theorem~$7'$]{Seidenberg1950}.
We can choose a nonzero section
$
    s\in H^0(\mathcal{X},\mathcal{N}_1^{\otimes m})$, 
such that $\div(s)$ is irreducible and normal. Write $\mathcal{Y}=\div(s)$, which is a normal projective variety of dimension $n-1$.
 Now we compute the interestion numbers.
\begin{equation*}
    m\left(\mathcal{M}^{n-r}\cdot
    \mathcal{N}_1\cdots\mathcal{N}_r\right)
    =
    \left(\mathcal{M}^{n-r}\cdot
    \mathcal{N}_2\cdots\mathcal{N}_r\cdot \mathcal{Y}\right)
    =
    \left((\mathcal{M}|_{\mathcal{Y}})^{n-r}\cdot
    \mathcal{N}_2|_{\mathcal{Y}}\cdots\mathcal{N}_r|_{\mathcal{Y}}\right),
\end{equation*}
and
\begin{equation*}
    m\left(\mathcal{M}^{n-r-1}\cdot\mathcal{L}\cdot
    \mathcal{N}_1\cdots\mathcal{N}_r\right)
    =
    \left(\mathcal{M}^{n-r-1}\cdot\mathcal{L}\cdot
    \mathcal{N}_2\cdots\mathcal{N}_r\cdot \mathcal{Y}\right)
    =
    \left((\mathcal{M}|_Y)^{n-r-1}\cdot\mathcal{L}|_{\mathcal{Y}}\cdot
    \mathcal{N}_2|_Y\cdots\mathcal{N}_r|_Y\right).
\end{equation*}
Then we apply the induction hypothesis for $r-1$ to $\mathcal{Y}$ and nef line bundles $\mathcal{L}|_{\mathcal{Y}}$, $\mathcal{M}|_{\mathcal{Y}}$,
$\mathcal{N}_2|_{\mathcal{Y}}$, $\dots$, $\mathcal{N}_r|_{\mathcal{Y}}$, we have
\begin{equation*}
    \mathcal{L}|_{\mathcal{Y}}\cdot \mathcal{N}_2|_{\mathcal{Y}}\cdots \mathcal{N}_r|_{\mathcal{Y}}\leq \mathcal{M}|_{\mathcal{Y}}\cdot \mathcal{N}_2|_{\mathcal{Y}}\cdots \mathcal{N}_r|_{\mathcal{Y}}.
\end{equation*}
For any nef line bundles $\mathcal{H}_1,\dots,\mathcal{H}_{n-r-1}$, their restrictions to $\mathcal{Y}$ are also nef. One gets
\begin{equation*}
\begin{aligned}
    \left(
        \mathcal{L}\cdot\mathcal{N}_1\cdots\mathcal{N}_r
        \cdot\mathcal{H}_1\cdots\mathcal{H}_{n-r-1}
    \right)
    &=
    \frac{1}{m}\left(
        \mathcal{L}|_Y\cdot
        \mathcal{N}_2|_Y\cdots\mathcal{N}_r|_Y
        \cdot\mathcal{H}_1|_Y\cdots\mathcal{H}_{n-r-1}|_Y
    \right)\\
    &\leq
    \frac{1}{m}\left(
        \mathcal{M}|_Y\cdot
        \mathcal{N}_2|_Y\cdots\mathcal{N}_r|_Y
        \cdot\mathcal{H}_1|_Y\cdots\mathcal{H}_{n-r-1}|_Y
    \right)\\
    &=
    \left(
        \mathcal{M}\cdot\mathcal{N}_1\cdots\mathcal{N}_r
        \cdot\mathcal{H}_1\cdots\mathcal{H}_{n-r-1}
    \right).
\end{aligned}
\end{equation*}
which completes the induction.
\end{proof}

\begin{corollary}\label{num corollary}
  Let $\mathcal{X}\to k$ be a normal projective variety with dimension $n$, $0\leq r\leq n-1$, $\mathcal{L}$, $\mathcal{M}$, $\mathcal{N}_1,\dots, N_r$ be nef line bundles, satisfies $\left(\mathcal{M}^{n-r}\cdot \mathcal{N}_1\cdots \mathcal{N}_r\right)>0$, then we have 
  \begin{equation*}
    \mathcal{L}\cdot \mathcal{N}_1\cdots \mathcal{N}_r\leq \left(n-r\right)\dfrac{\left(\mathcal{M}^{n-r-1}\cdot \mathcal{L}\cdot \mathcal{N}_1\cdots \mathcal{N}_r\right)}{\left(M^{n-r}\cdot N_1\cdots N_r\right)}\mathcal{M}\cdot \mathcal{N}_1\cdots N_r
  \end{equation*} 
\end{corollary}
\begin{proof}
We take 
$
    \lambda>
    (n-r)
    \dfrac{
        \left(\mathcal{M}^{n-r-1}\cdot\mathcal{L}\cdot
        \mathcal{N}_1\cdots\mathcal{N}_r\right)
    }{
        \left(\mathcal{M}^{n-r}\cdot
        \mathcal{N}_1\cdots\mathcal{N}_r\right)
    }.$ so $\lambda\mathcal{M}$
is nef.  Since
\begin{equation*}
\begin{aligned}
    &\left((\lambda\mathcal{M})^{n-r}\cdot
    \mathcal{N}_1\cdots\mathcal{N}_r\right)
    -(n-r)
    \left((\lambda\mathcal{M})^{n-r-1}\cdot
    \mathcal{L}\cdot\mathcal{N}_1\cdots\mathcal{N}_r\right)
    \\
    &\qquad=
    \lambda^{n-r-1}
    \left(
        \lambda\left(\mathcal{M}^{n-r}\cdot
        \mathcal{N}_1\cdots\mathcal{N}_r\right)
        -(n-r)
        \left(\mathcal{M}^{n-r-1}\cdot\mathcal{L}\cdot
        \mathcal{N}_1\cdots\mathcal{N}_r\right)
    \right)
    >0.
\end{aligned}
\end{equation*}
So we can apply Lemma \ref{num inequality} to
$\mathcal{L},\lambda\mathcal{M},
\mathcal{N}_1,\dots,\mathcal{N}_r$, one gets
\begin{equation*}
    \mathcal{L}\cdot\mathcal{N}_1\cdots\mathcal{N}_r
    \leq
    \lambda\mathcal{M}\cdot
    \mathcal{N}_1\cdots\mathcal{N}_r.
\end{equation*}
Letting $\lambda\to(n-r)
    \dfrac{
        \left(\mathcal{M}^{n-r-1}\cdot\mathcal{L}\cdot
        \mathcal{N}_1\cdots\mathcal{N}_r\right)
    }{
        \left(\mathcal{M}^{n-r}\cdot
        \mathcal{N}_1\cdots\mathcal{N}_r\right)
    }$, 
then we obtain the
desired inequality immediately
 \begin{equation*}
    \mathcal{L}\cdot \mathcal{N}_1\cdots \mathcal{N}_r\leq \left(n-r\right)\dfrac{\left(\mathcal{M}^{n-r-1}\cdot \mathcal{L}\cdot \mathcal{N}_1\cdots \mathcal{N}_r\right)}{\left(M^{n-r}\cdot \mathcal{N}_1\cdots \mathcal{N}_r\right)}\mathcal{M}\cdot \mathcal{N}_1\cdots \mathcal{N}_r.
  \end{equation*} 
\end{proof}

\begin{lemma}\label{higher codimensional siu inequality} Let $\mathcal{X}\to k$ be a normal projective variety with dimension $n$
    ,  $r \geq 0$ and $s \geq 1$ be integers such that $r+s \leq n$.
    Let $\mathcal{L}_1,\dots,\mathcal{L}_s,\mathcal{M},\mathcal{N}_1,\dots,N_r$ be nef
    line bundles , satisfies  $\left(\mathcal{M}^{n-r}\cdot \mathcal{N}_1\cdots \mathcal{N}_r\right)>0$. Then we have 
         \begin{equation*}
        \mathcal{L}_1\cdots \mathcal{L}_s \cdot \mathcal{N}_1\cdots \mathcal{N}_r\leq \left(n-r-s+1\right)^s \dfrac{\left(\mathcal{M}^{n-r-s}\cdot \mathcal{L}\cdots \mathcal{L}_s \cdot \mathcal{N}_1\cdots \mathcal{N}_r\right)}{\mathcal{M}^{n-r}\cdot \mathcal{N}_1\cdots \mathcal{N}_r} \mathcal{M}^s \cdot \mathcal{N}_1\cdots \mathcal{N}_r
    \end{equation*}

\end{lemma}
\begin{proof}
  By approximating the given nef line bundles by ample
line bundles, we may assume that all the line bundles are ample. Thus,  replacing $\mathcal{N}_1\cdots\mathcal{N}_r$ by $\mathcal{M}^{i-1}\cdot \mathcal{L}_{i+1}\cdots\mathcal{L}_s\cdot\mathcal{N}_1\cdots\mathcal{N}_r$ and $\mathcal{L}$ by $\mathcal{L}_i$ in the corollary \ref{num corollary} for $1\leq i\leq s$, and because $\left(\mathcal{M}^{n-r}\cdot \mathcal{N}_1\cdots \mathcal{N}_r\right)>0$, we obtained
\begin{equation*}
\begin{aligned}
    &\mathcal{M}^{i-1}\cdot
    \mathcal{L}_i\cdots\mathcal{L}_s\cdot
    \mathcal{N}_1\cdots\mathcal{N}_r
    \\
    &\qquad\leq
    (n-r-s+1)
    \dfrac{
        \left(\mathcal{M}^{n-r-s+i-1}\cdot
        \mathcal{L}_i\cdots\mathcal{L}_s\cdot
        \mathcal{N}_1\cdots\mathcal{N}_r\right)
    }{
        \left(\mathcal{M}^{n-r-s+i}\cdot
        \mathcal{L}_{i+1}\cdots\mathcal{L}_s\cdot
        \mathcal{N}_1\cdots\mathcal{N}_r\right)
    }
    \mathcal{M}^i\cdot
    \mathcal{L}_{i+1}\cdots\mathcal{L}_s\cdot
    \mathcal{N}_1\cdots\mathcal{N}_r.
\end{aligned}
\end{equation*}
Applying these inequalities  for each
$i=1,\dots,s$, we get
\begin{equation*}
    \begin{array}{ll}
        \mathcal{L}_{1}\cdots& \mathcal{L}_{s}\cdot \mathcal{N}_1\cdots \mathcal{N}_r\\
        &\leq \left(N-r-s+1\right)\dfrac{\left(\mathcal{M}^{n-r-s}\cdot \mathcal{L}_1 \cdots \mathcal{L}_s\cdot \mathcal{N}_1\cdots \mathcal{N}_r\right)}{\left(\mathcal{M}^{n-r-s+1}\cdot \mathcal{L}_{2}\cdots \mathcal{L}_s\cdot \mathcal{N}_1\cdots \mathcal{N}_r\right)} \mathcal{M}\cdot \mathcal{L}_{2}\cdots \mathcal{L}_{s}\cdot \mathcal{N}_1\cdots \mathcal{N}_r\\
        &\leq\left(N-r-s+1\right)^2\dfrac{\left(\mathcal{M}^{n-r-s}\cdot \mathcal{L}_1 \cdots \mathcal{L}_s\cdot \mathcal{N}_1\cdots \mathcal{N}_r\right)}{\left(\mathcal{M}^{n-r-s+2}\cdot \mathcal{L}_{3}\cdots \mathcal{L}_s\cdot \mathcal{N}_1\cdots \mathcal{N}_r\right)} \mathcal{M}^2\cdot \mathcal{L}_{3}\cdots \mathcal{L}_{s}\cdot \mathcal{N}_1\cdots \mathcal{N}_r\\
        &\leq\cdots\\
        &\leq\left(N-r-s+1\right)^s\dfrac{\left(\mathcal{M}^{n-r-s}\cdot \mathcal{L}_1 \cdots \mathcal{L}_s\cdot \mathcal{N}_1\cdots \mathcal{N}_r\right)}{\left(\mathcal{M}^{n-r}\cdot \mathcal{N}_1\cdots \mathcal{N}_r\right)} \mathcal{M}^{s}\cdot \mathcal{N}_1\cdots \mathcal{N}_r.
    \end{array}
\end{equation*}
which completes the proof.
\end{proof}

\begin{lemma}\label{key inequality}
    Let $\pi:\mathcal{X}\to B$ define as above with relative dimension $m$, $r,s$ be integers satisfy $r\geq 0$, $s>0$, and $r+s\leq m$.  Let $H$ be an ample line bundle on $B$, $\mathcal{M},\mathcal{N},\mathcal{L}$ be nef line bundles
    on $\mathcal{X}$ with generic fibers 
    $M=\mathcal{M}_K$, $N=\mathcal{N}_K$, and $L=\mathcal{L}_K$. Denote
    \begin{equation*}
        \mathcal{M}(t)
        =\mathcal{M}\otimes\pi^{*}H^{\otimes t}.
    \end{equation*}
   Assume 
    \begin{equation*}
        \left(L^s\cdot M^{m-r-s}\cdot N^r\right)>0,
        \qquad
        \left(M^{m-r}\cdot N^r\right)>0.
    \end{equation*}
    Then there exists a positive integer $t>0$, depending only
    on $m$, $s$, $\mathcal{M}$, $\mathcal{N}$, $\mathcal{L}$ and $H$.  Such that for any normal
    projective variety $\mathcal{Y}\to k$ and  any dominant generically
    finite morphism $g:\mathcal{Y}\to\mathcal{X}$ , we have
    \begin{equation*}
        g^{*}\mathcal{L}^{s}\cdot g^{*}\mathcal{N}^{r}
        \leq
        \left(m+d+1-r-s\right)^s
        \dfrac{
            \left(L^s\cdot M^{m-r-s}\cdot N^r\right)
        }{
            \left(M^{m-r}\cdot N^r\right)
        }
        \cdot
        g^{*}\left(
            \mathcal{M}\otimes\pi^{*}H^{\otimes t}
        \right)^s
        \cdot g^{*}\mathcal{N}^{r}.
    \end{equation*}
\end{lemma}
\begin{proof}
   Applying the Lemma~\ref{higher codimensional siu inequality}
    to the nef line bundles $g^{*}\mathcal{M}(t)$, $s$ copies
    $g^{*}\mathcal{L}$
    and $r$ copies of $g^{*}\mathcal{N}$, we obtain
  \begin{equation*}
        \begin{array}{ll}
            g^{*}\mathcal{L}^s \cdot g^{*}\mathcal{N}^r&\\
            \leq \left(m+d+1-r-s\right)^s \dfrac{\left(g^{*}\mathcal{L}^s\cdot g^{*}\mathcal{M}(t)^{m+d-r-s}\cdot g^{*}\mathcal{N}^r\right)}{\left(g^{*}\mathcal{M}(t)^{m+d-r}\cdot g^{*}\mathcal{N}^{r}\right)}\cdot g^{*}\mathcal{M}(t)^{s} \cdot g^{*}\mathcal{N}^r&\\
            = \left(m+d+1-r-s\right)^s\dfrac{\left(\mathcal{L}^{s}\cdot\mathcal{M}(t)^{m+d-r-s}\cdot\mathcal{N}^r\right)}{\left(\mathcal{M}(t)^{m+d-r}\cdot \mathcal{N}^r\right)}\cdot g^{*}\mathcal{M}(t)^{s} \cdot g^{*}\mathcal{N}^r.&\\
        \end{array}
    \end{equation*}
On the other hand, by the binomial expansion , we have
\begin{equation*}
    \begin{array}{ll}
       \dfrac{\left(\mathcal{L}^{s}\cdot\mathcal{M}(t)^{m+d-r-s}\cdot\mathcal{N}^r\right)}{\left(\mathcal{M}(t)^{m+d-r}\cdot \mathcal{N}^r\right)}
        =\dfrac{\sum\limits_{j=0}^{d}\binom{m+d-r-s}{j} t^{j}\left(\mathcal{L}^{s}\cdot\mathcal{M}^{m+d-r-s-j}\cdot\mathcal{N}^r\cdot (\pi^{*}H)^{j}\right)}{\sum\limits_{i=0}^{d}\binom{m+d-r}{i}t^{i}\left(\mathcal{M}^{m+d-r-i}\cdot \mathcal{N}^r\cdot \left(\pi^{*}H\right)^{i}\right)}.&\\
    \end{array}
\end{equation*}
Since the ratio of the leading coefficients is $\frac{\binom{m+d-r-s}{d}}{\binom{m+d-r}{d}}
    \frac{\left(\mathcal{L}^{s}\cdot\mathcal{M}^{m-r-s}\cdot\mathcal{N}^r\cdot (\pi^{*}H)^{d}\right)}{\left(\mathcal{M}^{m-r}\cdot \mathcal{N}^r\cdot \left(\pi^{*}H\right)^{d}\right)}$, and  $\frac{\binom{m+d-r-s}{d}}{\binom{m+d-r}{d}}<1$. 
We can take t sufficiently large such that
\begin{equation*}
    \dfrac{\left(\mathcal{L}^{s}\cdot\mathcal{M}(t)^{m+d-r-s}\cdot\mathcal{N}^r\right)}{\left(\mathcal{M}(t)^{m+d-r}\cdot \mathcal{N}^r\right)}
    \leq \frac{\left(\mathcal{L}^{s}\cdot\mathcal{M}^{m-r-s}\cdot\mathcal{N}^r\cdot (\pi^{*}H)^{d}\right)}{\left(\mathcal{M}^{m-r}\cdot \mathcal{N}^r\cdot \left(\pi^{*}H\right)^{d}\right)}=\dfrac{\left(L^s\cdot M^{m-r-s}\cdot N^r\right)}{ \left(M^{m-r}\cdot N^r\right)}.
\end{equation*}
Thus, we have
\begin{equation*}
    g^{*}\mathcal{L}^{s}\cdot g^{*}\mathcal{N}^{r}\leq \left(m+d+1-r-s\right)^s \dfrac{\left(L^s\cdot M^{m-r-s}\cdot N^r\right)}{ \left(M^{m-r}\cdot N^r\right)}\cdot  g^{*}\left(\mathcal{M}\otimes \pi^{*} H^{\otimes t}\right)^{s}\cdot g^{*}\mathcal{N}^{r}
\end{equation*}
\end{proof}

\begin{lemma}\label{inequality between line bundles and subvarieties}
  Let $\mathcal{X}/k$ be a normal projective variety of dimension $n$,
and let $1\leq r\leq n$ be an integer.
Let $\mathcal{L}_1,\dots,\mathcal{L}_r$ be nef line bundles on
$\mathcal{X}$, and suppose that there exist r nonzero sections
$s_1\in H^0(\mathcal{X},\mathcal{L}_1),\dots,s_r\in H^0(\mathcal{X},\mathcal{L}_r)$ for $1\leq i\leq r$.
Let $\mathcal{V}$ be an irreducible component of
$\ddiv(s_1)\cap\cdots\cap\ddiv(s_r)$ of codimension $r$ in
$\mathcal{X}$.
Then, for any nef line bundles
$\mathcal{H}_1,\dots,\mathcal{H}_{n-r}$ on $\mathcal{X}$, we have
\begin{equation*}
    \left(
    \mathcal{H}_1\cdots\mathcal{H}_{n-r}
    \cdot\mathcal{L}_1\cdots\mathcal{L}_r
    \right)
    \geq
    \left(
    \mathcal{H}_1\cdots\mathcal{H}_{n-r}\cdot\mathcal{V}
    \right).
\end{equation*}
\end{lemma}
\begin{proof}
  We prove the lemma by induction on $r$.
For $r=1$, $\mathcal{V}$ is an irreducible component of
$\ddiv(s_1)$. Since $[\ddiv(s_1)]-[\mathcal{V}]$ is an effective
cycle, we have
\begin{equation*}
    \left(
    \mathcal{H}_1\cdots\mathcal{H}_{n-1}\cdot\mathcal{L}_1
    \right)
    =
    \left(
    \mathcal{H}_1\cdots\mathcal{H}_{n-1}\cdot\ddiv(s_1)
    \right)
    \geq
    \left(
    \mathcal{H}_1\cdots\mathcal{H}_{n-1}\cdot\mathcal{V}
    \right).
\end{equation*}
Suppose that the lemma holds for $r=s-1$, and consider the case
$r=s$. Fix any nef line bundles
$\mathcal{H}_1,\dots,\mathcal{H}_{n-s}$. We need to show that
\begin{equation*}
    \left(
    \mathcal{H}_1\cdots\mathcal{H}_{n-s}
    \cdot\mathcal{L}_1\cdots\mathcal{L}_s
    \right)
    \geq
    \left(
    \mathcal{H}_1\cdots\mathcal{H}_{n-s}\cdot\mathcal{V}
    \right).
\end{equation*}
Now assume that $r\geq 2$ and that the lemma holds for $r-1$. As in the proof of \cite[Lemma~3.3]{JSXZ2024},   There exists an irreducible component $\mathcal{W}$ of
$\bigcap\limits_{i=1}^{r-1}\ddiv(s_i)$ of codimension $r-1$
in $\mathcal{X}$ such that
$\mathcal{W}\nsubseteq\ddiv(s_r)$ and $\mathcal{V}$ is an
irreducible component of $\mathcal{W}\cap\ddiv(s_r)$.
Then by the induction hypothesis, applied with the nef line bundles
$\mathcal{H}_1,\dots,\mathcal{H}_{n-r},\mathcal{L}_r$, we have
\begin{equation*}
    \begin{array}{rl}
        \left(
        \mathcal{H}_1\cdots\mathcal{H}_{n-r}
        \cdot\mathcal{L}_1\cdots\mathcal{L}_r
        \right)
        &\geq
        \left(
        \mathcal{H}_1\cdots\mathcal{H}_{n-r}
        \cdot\mathcal{L}_r\cdot\mathcal{W}
        \right)\\
        &=
        \left(
        \mathcal{H}_1\cdots\mathcal{H}_{n-r}
        \cdot\mathcal{W}\cdot\ddiv(s_r)
        \right)\\
        &\geq
        \left(
        \mathcal{H}_1\cdots\mathcal{H}_{n-r}\cdot\mathcal{V}
        \right).
    \end{array}
\end{equation*}
which completes the proof.
\end{proof}

\subsection{Arithmetic dynamical degrees}
Let $K$ be a finitely generated extension of $k$
with $\trdeg_k(K)=d\geq 1$.
Let $X$ be a normal projective variety over $K$
with dimension $n$, and 
$f:X\dashrightarrow X$ be a dominant rational self-map.
Define
\begin{equation*}
  X_f(\overline{K})
    :=
    \left\{
        x\in X(\overline{K})
        \mid f^n(x)\notin I(f) \text{ for all } n\geq 0
    \right\},
\end{equation*}
where $I(f^n)$ denotes the indeterminacy locus of $f^n$.

Let $\mathcal{B}=(B;H,\dots,H)$ be a polarization of $K$,
where $H$ is an ample line bundle on $B$ and occurs
$d-1$ times.
Let $L$ be an ample line bundle on $X$, and consider
a normal projective $B$-model $(\mathcal{X},\mathcal{L})$
of $(X,L)$, that is, $\mathcal{X}$ is a $B$-model of  $X$, and $\mathcal{L}$ is an ample line bundle on $\mathcal{X}$  such that $\mathcal{L}_K=L$.

The dominant rational self-map $f:X\dashrightarrow X$
induces a dominant rational self-map
$F:\mathcal{X}\dashrightarrow\mathcal{X}$ over $B$.
Taking the normalization of the closure of its graph,
we obtain the following diagram:
\begin{center}
   \begin{tikzcd}
 & \Gamma_{F} \arrow{dl}[swap]{\pi_{1}} \arrow{dr}{\pi_2}\\
\mathcal{X} \arrow[dashed]{rr}{F} && \mathcal{X} ,
\end{tikzcd} 
\end{center} 
where $\Gamma_F$ is the normalization of the closure
of the graph of $F$ in $\mathcal{X}\times_B\mathcal{X}$,
and $\pi_1,\pi_2$ are the projections to the first and
second factors, respectively.
The morphism $\pi_1$ is birational, and $\pi_2$ is
dominant and generically finite.
We also use $\pi_1$ and $\pi_2$ to denote the induced
morphisms on the generic fibers over $B$.

\begin{definition}
    For an integer $1\leq k\leq n$, we define
    \begin{equation*}
        \widehat{\deg}_{k,\mathcal{L}}(F)
        :=
        \left(
            \pi_1^{*}\mathcal{L}^{n-k+1}
            \cdot\pi_2^{*}\mathcal{L}^{k}
            \cdot\pi_1^{*}\pi^{*}H^{d-1}
        \right).
    \end{equation*}
\end{definition}

By the choice of $\mathcal{L}$ and $H$,
$\widehat{\deg}_{k,\mathcal{L}}(F)$ is positive.
We next prove that these arithmetic degrees satisfy
the following submultiplicative inequality.
\begin{theorem}\label{multi of degree}
    With the notation above, there exists a constant
    $C>0$ such that for any dominant rational self-maps
    $F:\mathcal{X}\dashrightarrow\mathcal{X}$ and
    $G:\mathcal{X}\dashrightarrow\mathcal{X}$ over $B$, we have
    \begin{equation*}
        \widehat{\deg}_{k,\mathcal{L}}(F\circ G)
        \leq C\cdot
        \widehat{\deg}_{k,\mathcal{L}}(F)\cdot
        \widehat{\deg}_{k,\mathcal{L}}(G).
    \end{equation*}
\end{theorem}
\begin{proof}
Consider the following diagram:
\begin{equation}\label{resolution of compose map}
    \begin{tikzcd}
        & & \widetilde{\Gamma}
        \arrow{ld}[swap]{\pi_3}
        \arrow{rd}{\widetilde{\pi}_3} & & \\
        & \Gamma_G
        \arrow{ld}[swap]{\widetilde{\pi}_1}
        \arrow{rd}{\widetilde{\pi}_2}
        & & \Gamma_F
        \arrow{ld}[swap]{\pi_1}
        \arrow{rd}{\pi_2} & \\
        \mathcal{X}
        \arrow{rrd}[swap]{\pi}
        \arrow[dashed]{rr}[swap]{G}
        & & \mathcal{X}
        \arrow{d}{\pi}
        \arrow[dashed]{rr}[swap]{F}
        & & \mathcal{X}
        \arrow{lld}{\pi} \\
        & & B & &
    \end{tikzcd}
\end{equation}
Here $\widetilde{\Gamma}$ is the normalization of the graph of $\pi_1^{-1}\circ\tilde{\pi}_2$ in $\Gamma_G\times_B\Gamma_F$. We denote $\pi_3$ and $\widetilde{\pi}_3$ denote the morphisms to
the first and second coordinates, respectively.
Then $\pi_3$ is birational, and
$\widetilde{\pi}_3$ is dominant and generically finite. 

Write $   \widetilde{\pi}
    :=
    \pi\circ\widetilde{\pi}_1\circ\pi_3
    =
    \pi\circ\pi_1\circ\widetilde{\pi}_3.$ By the projection formula, we have
\begin{equation*}
    \left(
        \pi_3^{*}\widetilde{\pi}_2^{*}\mathcal{L}^{n+1}
        \cdot\widetilde{\pi}^{*}H^{d-1}
    \right)
    =
    \deg(\widetilde{\pi}_2)
    \left(
        \mathcal{L}^{n+1}\cdot\pi^{*}H^{d-1}
    \right)>0.
\end{equation*}
Applying Lemma \ref{higher codimensional siu inequality}
on $\widetilde{\Gamma}$ with $r=d-1$ to
$\widetilde{\pi}_3^{*}\pi_2^{*}\mathcal{L}$,
$\pi_3^{*}\widetilde{\pi}_2^{*}\mathcal{L}$,
and $d-1$ copies of $\widetilde{\pi}^{*}H$, we obtain
\begin{equation*}
        \widetilde{\pi}_3^{*}\pi_2^{*}\mathcal{L}^k
        \cdot\widetilde{\pi}^{*}H^{d-1}\leq
        (n+2-k)^k
        \dfrac{
            \left(
                \pi_3^{*}\widetilde{\pi}_2^{*}\mathcal{L}^{n+1-k}
                \cdot\widetilde{\pi}_3^{*}\pi_2^{*}\mathcal{L}^k
                \cdot\widetilde{\pi}^{*}H^{d-1}
            \right)
        }{
            \left(
                \pi_3^{*}\widetilde{\pi}_2^{*}\mathcal{L}^{n+1}
                \cdot\widetilde{\pi}^{*}H^{d-1}
            \right)
        }
        \cdot
        \pi_3^{*}\widetilde{\pi}_2^{*}\mathcal{L}^k
        \cdot\widetilde{\pi}^{*}H^{d-1}.
\end{equation*}
Intersecting both sides of the inequality above with with
$\pi_3^{*}\widetilde{\pi}_1^{*}\mathcal{L}^{n+1-k}$, one gets
\begin{equation*}
   \widehat{\deg}_{k,\mathcal{L}}(F\circ G)\leq
        (n+2-k)^k
        \dfrac{
            \left(
                \pi_3^{*}\widetilde{\pi}_2^{*}\mathcal{L}^{n+1-k}
                \cdot\widetilde{\pi}_3^{*}\pi_2^{*}\mathcal{L}^k
                \cdot\widetilde{\pi}^{*}H^{d-1}
            \right)
        }{
            \left(
                \pi_3^{*}\widetilde{\pi}_2^{*}\mathcal{L}^{n+1}
                \cdot\widetilde{\pi}^{*}H^{d-1}
            \right)
        }
        \cdot\widehat{\deg}_{k,\mathcal{L}}(G).
\end{equation*}
On the other hand, since
$\pi_1\circ\widetilde{\pi}_3
=\widetilde{\pi}_2\circ\pi_3$,
we obtain
\begin{equation*}
    \begin{aligned}
        \dfrac{
            \left(
                \pi_3^{*}\widetilde{\pi}_2^{*}\mathcal{L}^{n+1-k}
                \cdot\widetilde{\pi}_3^{*}\pi_2^{*}\mathcal{L}^k
                \cdot\widetilde{\pi}^{*}H^{d-1}
            \right)
        }{
            \left(
                \pi_3^{*}\widetilde{\pi}_2^{*}\mathcal{L}^{n+1}
                \cdot\widetilde{\pi}^{*}H^{d-1}
            \right)
        }&=
        \dfrac{
            \left(
                \widetilde{\pi}_3^{*}\pi_1^{*}\mathcal{L}^{n+1-k}
                \cdot\widetilde{\pi}_3^{*}\pi_2^{*}\mathcal{L}^k
                \cdot\widetilde{\pi}_3^{*}\pi_1^{*}\pi^{*}H^{d-1}
            \right)
        }{
            \left(
                \widetilde{\pi}_3^{*}\pi_1^{*}\mathcal{L}^{n+1}
                \cdot\widetilde{\pi}_3^{*}\pi_1^{*}\pi^{*}H^{d-1}
            \right)
        }\\
        &=
        \dfrac{
            \deg(\widetilde{\pi}_3)
            \left(
                \pi_1^{*}\mathcal{L}^{n+1-k}
                \cdot\pi_2^{*}\mathcal{L}^k
                \cdot\pi_1^{*}\pi^{*}H^{d-1}
            \right)
        }{
            \deg(\widetilde{\pi}_3)
            \left(
                \pi_1^{*}\mathcal{L}^{n+1}
                \cdot\pi_1^{*}\pi^{*}H^{d-1}
            \right)
        }\\
        &=
        \dfrac{
            \widehat{\deg}_{k,\mathcal{L}}(F)
        }{
            \left(\mathcal{L}^{n+1}\cdot\pi^{*}H^{d-1}\right)
        }.
    \end{aligned}
\end{equation*}
Set
$C:=
    \dfrac{
        (n+2-k)^k
    }{
        \left(\mathcal{L}^{n+1}\cdot\pi^{*}H^{d-1}\right)
    }.$
Thus, there is a constant independent of $F$ and $G$, such that
\begin{equation*}
    \widehat{\deg}_{k,\mathcal{L}}(F\circ G)
    \leq C\cdot
    \widehat{\deg}_{k,\mathcal{L}}(F)\cdot
    \widehat{\deg}_{k,\mathcal{L}}(G).
\end{equation*}
\end{proof}

To prove the existence of arithmetic dynamical degrees,
we recall the following well-known lemma of Fekete.
\begin{lemma}[Fekete] Let $\{a_n\}_{n\geq 0}$ be any submultiplicative sequence of positive real numbers, i.e. for any $m,n\geq 0$
\begin{equation*}
  a_{m+n}\leq a_m\cdot a_n.
\end{equation*}
Then the limit $\lim\limits_{n\to \infty} a_n^{1/n}$ exists. 
\end{lemma}
\begin{proof}
Set $a_0:=1$ and
\begin{equation*}
    \ell:=\inf_{n\geq 1}a_n^{1/n}.
\end{equation*}
Then
\begin{equation*}
    \ell\leq\liminf_{n\to\infty}a_n^{1/n}.
\end{equation*}

Fix an integer $m\geq 1$.
For each $n\geq 1$, write $n=qm+r$, where
$q\geq 0$ and $0\leq r<m$.
By submultiplicativity, we have
\begin{equation*}
    a_n\leq a_m^q a_r\leq a_m^q C_m,
\end{equation*}
where
\begin{equation*}
    C_m:=\max\{a_0,a_1,\dots,a_{m-1}\}.
\end{equation*}
Taking $n$-th roots, we obtain
\begin{equation*}
    a_n^{1/n}\leq a_m^{q/n}C_m^{1/n}.
\end{equation*}
Since $q/n\to 1/m$ as $n\to\infty$, it follows that
\begin{equation*}
    \limsup_{n\to\infty}a_n^{1/n}\leq a_m^{1/m}.
\end{equation*}
Taking the infimum over all $m\geq 1$, we have
\begin{equation*}
    \limsup_{n\to\infty}a_n^{1/n}\leq\ell.
\end{equation*}
Therefore, the limit
\begin{equation*}
    \lim_{n\to\infty}a_n^{1/n}=\ell.
\end{equation*}
exists.
\end{proof}

By Theorem \ref{multi of degree}, the sequence
$\{C\widehat{\deg}_{k,\mathcal{L}}(F^n)\}_{n\geq 0}$
is submultiplicative.
Applying Fekete's lemma and noting that $C^{1/n}\to 1$,
we see that the limit
\begin{equation*}
    \lim_{n\to\infty}
    \widehat{\deg}_{k,\mathcal{L}}(F^n)^{1/n}
\end{equation*}
exists.
Therefore, we can define the arithmetic dynamical degrees
as follows.

\begin{definition}
 For an integer $1 \leq k\leq n $ , we define the $k$-th arithmetic dynamical degree
  \begin{equation*}
        \widehat{\lambda}_k\left(f\right):=\lim\limits_{m\to \infty}\widehat{\deg}_{k,\mathcal{L}}\left(F^m\right)^{1/m}.
    \end{equation*}
\end{definition}

\begin{theorem}\label{relation between arithmetic dynamical degree and dynamical degree}
 Let $f:X\dashrightarrow X$ be a dominant rational self-map,
and let $1\leq k\leq n$ be an integer.
With the notation defined as above, we have
  \begin{equation*}
  \widehat{\lambda}_k\left(f\right)=\max\{\lambda_k(f), \lambda_{k-1}(f)\}.
  \end{equation*}
\end{theorem}
\begin{proof}
  We first prove
\begin{equation*}
    \widehat{\lambda}_{k}(f)\geq \max\left\{\lambda_k(f),\lambda_{k-1}(f)\right\}.
\end{equation*}
Applying Lemma \ref{higher codimensional siu inequality} to
$\pi_1^{*}\mathcal{L}$ and
$\pi_1^{*}\pi^{*}H$, we have
\begin{equation*}
    \pi_1^{*}\pi^{*}H\leq\left(n+d\right)\dfrac{
\left(\pi_1^{*}\pi^{*}H\cdot\pi_1^{*}\mathcal{L}^{n+d-1}\right)}{\left(\pi_1^{*}\mathcal{L}^{n+d}\right)}\cdot\pi_1^{*}\mathcal{L}.
\end{equation*}
We intersect both sides of the inequality with
$\pi_1^{*}\mathcal{L}^{n-k}
\cdot\pi_2^{*}\mathcal{L}^{k}
\cdot\pi_1^{*}\pi^{*}H^{d-1}$,
we obtain
\begin{equation*}
    \begin{array}{rl}
          \widehat{\deg}_{k,\mathcal{L}}(F)=&\left(\pi_1^{*}\mathcal{L}^{n-k+1}\cdot \pi_2^{*}\mathcal{L}^{k}\cdot \pi_1^{*}\pi^{*}H^{d-1}\right)\\
          &\geq \dfrac{\left(\pi_1^{*}\mathcal{L}^{n+d}\right)}{\left(n+d\right)\left(\pi_1^{*}\pi^{*}H\cdot \pi_1^{*}\mathcal{L}^{n+d-1}\right)}\cdot \left(\pi_1^{*}\mathcal{L}^{n-k}\cdot \pi_2^{*}\mathcal{L}^{k}\cdot \pi_1^{*}\pi^{*}H^{d}\right).
    \end{array}
\end{equation*}
On the other hand, the projection formula gives
\begin{equation*}
    \left(\pi_1^{*}\mathcal{L}^{n-k}\cdot \pi_2^{*}\mathcal{L}^{k}\cdot \pi_1^{*}\pi^{*}H^{d}\right)=\left(H^d\right)\left( \pi_2^{*}L^k\cdot \pi_1^{*}L^{n-k}\right).
\end{equation*}
Set
   $ C_0:=
    \dfrac{
        \left(\mathcal{L}^{n+d}\right)(H^d)
    }{
        (n+d)\left(\pi^{*}H\cdot\mathcal{L}^{n+d-1}\right)
    }>0.$
Then
\begin{equation*}
    \widehat{\deg}_{k,\mathcal{L}}(F)
    \geq
    C_0\deg_{k,L}(F_K,\mathcal{X}_K/B_K).
\end{equation*}
Similarly, we consider the second projection.
Since $\pi\circ\pi_1=\pi\circ\pi_2$, we have
\begin{equation*}
    \pi_1^{*}\pi^{*}H=\pi_2^{*}\pi^{*}H.
\end{equation*}
So by the projection formula, we have
\begin{equation*}
        \dfrac{
            \left(
                \pi_1^{*}\pi^{*}H
                \cdot\pi_2^{*}\mathcal{L}^{n+d-1}
            \right)
        }{
            \left(\pi_2^{*}\mathcal{L}^{n+d}\right)
        }
        =
        \dfrac{
            \deg(\pi_2)
            \left(\pi^{*}H\cdot\mathcal{L}^{n+d-1}\right)
        }{
            \deg(\pi_2)\left(\mathcal{L}^{n+d}\right)
        }
        =
        \dfrac{
            \left(\pi^{*}H\cdot\mathcal{L}^{n+d-1}\right)
        }{
            \left(\mathcal{L}^{n+d}\right)
        }.
\end{equation*}
Applying Lemma \ref{higher codimensional siu inequality} to $\pi_1^{*}\pi^{*}H$ and $ \pi_2^{*}\mathcal{L}$, we obtain
\begin{equation*}
    \pi_1^{*}\pi^{*}H
    \leq
    (n+d)
    \dfrac{
        \left(\pi^{*}H\cdot\mathcal{L}^{n+d-1}\right)
    }{
        \left(\mathcal{L}^{n+d}\right)
    }
    \cdot\pi_2^{*}\mathcal{L}.
\end{equation*}
Intersecting both sides of this inequality with
$\pi_1^{*}\mathcal{L}^{n-k+1}
\cdot\pi_2^{*}\mathcal{L}^{k-1}
\cdot\pi_1^{*}\pi^{*}H^{d-1}$, we have
\begin{equation*}
  \begin{array}{rl}
     \widehat{\deg}_{k,\mathcal{L}}(F)&=\left(\pi_1^{*}\mathcal{L}^{n-k+1}\cdot \pi_2^{*}\mathcal{L}^{k}\cdot \pi_1^{*}\pi^{*}H^{d-1}\right)\\
        &\geq
        \dfrac{
            \left(\mathcal{L}^{n+d}\right)
        }{
            (n+d)
            \left(\pi^{*}H\cdot\mathcal{L}^{n+d-1}\right)
        }
        \left(
            \pi_1^{*}\mathcal{L}^{n-k+1}
            \cdot\pi_2^{*}\mathcal{L}^{k-1}
            \cdot\pi_1^{*}\pi^{*}H^d
        \right)\\
        &=
        C_0\deg_{k-1,L}(F_K,\mathcal{X}_K/B_K).
  \end{array}
\end{equation*}
The constant $C_0$ is independent of $F$. Now we
replace $F$ by $F^s$, and taking $s$-th roots,
letting $s\to\infty$, we obtain
\begin{equation*}
    \widehat{\lambda}_k(f)\geq\max\{\lambda_k(f),\lambda_{k-1}(f)\}.
\end{equation*}

Now we consider the converse inequality. We again
use the diagram.
\begin{equation*}
     \begin{tikzcd}
		&                                                   & \widetilde{\Gamma} \arrow{ld}[swap]{\pi_3} \arrow{rd}{\widetilde{\pi}_3}   &                                                   &                    \\
		& \Gamma_G \arrow{ld}[swap]{\widetilde{\pi}_1} \arrow{rd}{\widetilde{\pi}_2} &                                           & \Gamma_F \arrow{ld}[swap]{\pi_1} \arrow{rd}{\pi_2} &                    \\
		\mathcal{X} \arrow{rrd}[swap]{\pi} \arrow[dashed]{rr}[swap]{G} &                                                   & \mathcal{X} \arrow{d}{\pi} \arrow[dashed]{rr}[swap]{F} &                                                   & \mathcal{X} \arrow{lld}{\pi} \\
		&                                                   & B                              &                                                   &                   
	\end{tikzcd}   
\end{equation*} 
Apply Lemma \ref{key inequality} with $r=0$ to
$\Gamma_F\to B$, the nef line bundles
$\pi_2^{*}\mathcal{L}$ and $\pi_1^{*}\mathcal{L}$,
and the morphism $\widetilde{\pi}_3$, we obtain a
positive integer $t=t(F)$, independent of $G$, such that
\begin{equation*}
\widetilde{\pi}_3^{*}\pi_2^{*}\mathcal{L}^k\leq(n+d+1-k)^k\dfrac{\left(\pi_2^{*}L^k\cdot\pi_1^{*}L^{n-k}\right)}{(L^n)}\cdot\widetilde{\pi}_3^{*}\pi_1^{*} \left(\mathcal{L}\otimes\pi^{*}H^{\otimes t}\right)^k.
\end{equation*}
By intersecting the inequality above with $\pi_3^{*}\widetilde{\pi}_1^{*}\mathcal{L}^{n+1-k}\cdot \widetilde{\pi}^{*}H^{d-1}$, we have
\begin{equation*}
    \begin{aligned}
        \widehat{\deg}_{k,\mathcal{L}}(F\circ G)
        &\leq
        \left(
            \pi_3^{*}\widetilde{\pi}_1^{*}\mathcal{L}^{n+1-k}
            \cdot\widetilde{\pi}_3^{*}\pi_2^{*}\mathcal{L}^k
            \cdot\widetilde{\pi}^{*}H^{d-1}
        \right)\\
        &\leq(n+d+1-k)^k\dfrac{\left(\pi_2^{*}L^k\cdot\pi_1^{*}L^{n-k}\right)}{(L^n)}
        \left(
            \pi_3^{*}\widetilde{\pi}_1^{*}\mathcal{L}^{n+1-k}
            \cdot\widetilde{\pi}_3^{*}\pi_1^{*}
            \left(\mathcal{L}\otimes\pi^{*}H^{\otimes t}\right)^k
            \cdot\widetilde{\pi}^{*}H^{d-1}
        \right)\\
        &=C_1\deg_{k,L}(f)
        \left(
            \pi_3^{*}\widetilde{\pi}_1^{*}\mathcal{L}^{n+1-k}
            \cdot\widetilde{\pi}_3^{*}\pi_1^{*}
            \left(\mathcal{L}\otimes\pi^{*}H^{\otimes t}\right)^k
            \cdot\widetilde{\pi}^{*}H^{d-1}
        \right),
    \end{aligned}
\end{equation*}
where $C_1:=(n+d+1-k)^k/(L^n)$.
By the binomial expansion and the projection formula,
all terms containing more than $d$ factors of
$\widetilde{\pi}^{*}H$ have zero intersection number.
Thus
\begin{equation*}
    \begin{aligned}
        &\widehat{\deg}_{k,\mathcal{L}}(F\circ G)\leq C_1\deg_{k,L}(f)\left(
            \pi_3^{*}\widetilde{\pi}_1^{*}\mathcal{L}^{n+1-k}
            \cdot\widetilde{\pi}_3^{*}\pi_1^{*}
            \left(\mathcal{L}\otimes\pi^{*}H^{\otimes t}\right)^k
            \cdot\widetilde{\pi}^{*}H^{d-1}
        \right)\\
        &\quad\leq
        C_1\deg_{k,L}(f)
        \Bigg[
        \left(
            \pi_3^{*}\widetilde{\pi}_1^{*}\mathcal{L}^{n+1-k}
            \cdot\widetilde{\pi}_3^{*}\pi_1^{*}\mathcal{L}^k
            \cdot\widetilde{\pi}^{*}H^{d-1}
        \right)
        +kt
        \left(
            \pi_3^{*}\widetilde{\pi}_1^{*}\mathcal{L}^{n+1-k}
            \cdot\widetilde{\pi}_3^{*}\pi_1^{*}\mathcal{L}^{k-1}
            \cdot\widetilde{\pi}^{*}H^d
        \right)
        \Bigg].
    \end{aligned}
\end{equation*}
Consider the first term, since
$\widetilde{\pi}_2\circ\pi_3
=\pi_1\circ\widetilde{\pi}_3$
and $\pi_3$ is birational, the projection formula gives
\begin{equation*}
    \left(
        \pi_3^{*}\widetilde{\pi}_1^{*}\mathcal{L}^{n+1-k}
        \cdot\widetilde{\pi}_3^{*}\pi_1^{*}\mathcal{L}^k
        \cdot\widetilde{\pi}^{*}H^{d-1}
    \right)
    =
    \widehat{\deg}_{k,\mathcal{L}}(G).
\end{equation*}
On the other hand,
\begin{equation*}
        \left(
            \pi_3^{*}\widetilde{\pi}_1^{*}\mathcal{L}^{n+1-k}
            \cdot\widetilde{\pi}_3^{*}\pi_1^{*}\mathcal{L}^{k-1}
            \cdot\widetilde{\pi}^{*}H^d
        \right)=
        (H^d)
        \left(
            \widetilde{\pi}_1^{*}L^{n+1-k}
            \cdot\widetilde{\pi}_2^{*}L^{k-1}
        \right)=
        (H^d)\deg_{k-1,L}(g).
\end{equation*}
Set $C_2:=C_1(H^d)$, we obtain
\begin{equation*}
    \widehat{\deg}_{k,\mathcal{L}}(F\circ G)
    \leq
    C_1\deg_{k,L}(f)\widehat{\deg}_{k,\mathcal{L}}(G)
    +
    ktC_2\deg_{k,L}(f)\deg_{k-1,L}(g).
\end{equation*}
According to the submultiplicity formula of  degrees in \cite[Theorem~1]{Dang20}
, we have
\begin{equation*}
    \deg_{k-1,L}(f\circ g)
    \leq
    C_3\deg_{k-1,L}(f)\deg_{k-1,L}(g),
\end{equation*}
where the constant $C_3$ is independent of the choice of $f$ and $g$.
We now consider the following vector and matrix:
\begin{equation*}
    U(F):=
    \begin{bmatrix}
        \widehat{\deg}_{k,\mathcal{L}}(F)\\
        \deg_{k-1,L}(f)
    \end{bmatrix},
    \quad
    M(f):=
    \begin{bmatrix}
        C_1\deg_{k,L}(f)
        &
        ktC_2\deg_{k,L}(f)\\
        0
        &
        C_3\deg_{k-1,L}(f)
    \end{bmatrix}.
\end{equation*}
aking $G=F^{s-1}$, we have
\begin{equation*}
    U(F^s)\leq M(f)U(F^{s-1}).
\end{equation*}
Since all entries of $M(f)$ are nonnegative,
repeated application gives
\begin{equation*}
    U(F^s)\leq M(f)^{s-1}U(F).
\end{equation*}
For a vector $v=(v_1,v_2)^{\top}$ and a matrix
$M=(m_{ij})$, write
\begin{equation*}
    \|v\|_{\infty}:=\max\{|v_1|,|v_2|\},
    \qquad
    \|M\|_{\max}:=\max_{i,j}|m_{ij}|.
\end{equation*}
Then
\begin{equation*}
        \widehat{\deg}_{k,\mathcal{L}}(F^s)\leq
        \|U(F^s)\|_{\infty}\leq
        \|M(f)^{s-1}U(F)\|_{\infty}\leq
        2\|M(f)^{s-1}\|_{\max}\|U(F)\|_{\infty}.
\end{equation*}
Now we fix   $r$, and replacing $F$ by $F^r$, taking $rs$-th roots.
we obtain
\begin{equation*}
    \widehat{\deg}_{k,\mathcal{L}}(F^{sr})^{1/sr}
    \leq
    2^{1/sr}
    \|M(f^r)^{s-1}\|_{\max}^{1/sr}
    \|U(F^r)\|_{\infty}^{1/sr}.
\end{equation*}
Here the integer $t$ is independent of $s$.
So if we let $s\to\infty$, the spectral radius formula gives
\begin{equation*}
    \widehat{\lambda}_k(f)
    \leq
    \rho(M(f^r))^{1/r}.
\end{equation*}
Since $M(f^r)$ is upper triangular, its eigenvalues are
$C_1\deg_{k,L}(f^r)$ and
$C_3\deg_{k-1,L}(f^r)$.
So we have
\begin{equation*}
    \widehat{\lambda}_k(f)
    \leq
    \max\left\{
        C_1\deg_{k,L}(f^r),
        C_3\deg_{k-1,L}(f^r)
    \right\}^{1/r}.
\end{equation*}
Finally, taking the limit as $r\to\infty$ yields the desired inequality
\begin{equation*}
    \widehat{\lambda}_k(f)
    \leq
    \max\{\lambda_k(f),\lambda_{k-1}(f)\}.
\end{equation*}
Hence, we have the assertion
\end{proof}

The following corollary follows immediately from the theorem above.

\begin{corollary}
    Under the assumptions above, the arithmetic dynamical degree
    $\widehat{\lambda}_k(f)$ is independent of the choice of
    the polarization $\mathcal{B}$, the model $\mathcal{X}$,
    and the ample line bundle $\mathcal{L}$ on $\mathcal{X}$.
\end{corollary}
\subsection{Kawaguchi--Silverman conjecture}
 In this section, we turns to study the Kawaguchi--Silverman
conjecture over finitely generated fields of positive characteristic.

Let $k$ be a finite field and  $K/k$ be a finitely
generated field with $\trdeg_k(K)=d\geq 1$.
Let $X$ be a normal projective variety over defined over $K$, let $L$ be an ample line bundle on $X$,
and let $\mathcal{B}$ be a big and nef polarization of $K$.
Choose a model $(\mathcal{X},\mathcal{L})$ of $(X,L)$, and let
$h^{\mathcal{B}}_{(\mathcal{X},\mathcal{L})}$ be the associated
Moriwaki height. We define the arithmetic degrees of points
as follows.

\begin{definition}[Arithmetic degree]\label{Arithmetic degree}
  Let $f:X\dashrightarrow X$ be a dominant rational self-map. Keep the notation defined as above , Set $h^{+}_{(\mathcal{X},\mathcal{L})}
    :=
    \max\left\{
    h^{\mathcal{B}}_{(\mathcal{X},\mathcal{L})},1
    \right\}.$
For any $x\in X_f(\overline{K})$, we define the upper and lower
arithmetic degrees of $x$ under $f$ by
\begin{equation*}
    \begin{aligned}
        \overline{\alpha}_f(x)
        &:=
        \limsup\limits_{n\to\infty}
        \left(
        h^{+}_{(\mathcal{X},\mathcal{L})}(f^n(x))
        \right)^{1/n},\\
        \underline{\alpha}_f(x)
        &:=
        \liminf\limits_{n\to\infty}
        \left(
        h^{+}_{(\mathcal{X},\mathcal{L})}(f^n(x))
        \right)^{1/n}.
    \end{aligned}
\end{equation*}
If these two limits are equal, we call 
\begin{equation*}
    \alpha_f(x)
    :=
    \lim\limits_{n\to\infty}
    \left(
    h^{+}_{(\mathcal{X},\mathcal{L})}(f^n(x))
    \right)^{1/n}
\end{equation*}
the arithmetic degree of $x$ under $f$.
\end{definition}
\begin{remark}
By Theorem  \ref{comparable of  height under field extension},
the corresponding height functions are comparable under field
extensions and changes of polarizations and models.
Consequently, $\overline{\alpha}_f(x)$ and
$\underline{\alpha}_f(x)$ are independent of the choice of
the polarization $\mathcal{B}$ and the model
$(\mathcal{X},\mathcal{L})$.
\end{remark}

Now, We can establish th fundamental inequality between arithmetic degrees and dynamical degrees.

\begin{theorem}\label{fundemental inequality over fg case}
Let $K$ be a finitely generated field over $\mathbb{F}_q$,
let $X$ be a normal projective variety over $K$,
and let $f:X\dashrightarrow X$ be a dominant rational self-map.
Then, we have
\begin{equation*}
    \overline{\alpha}_f(x)\leq \widehat{\lambda}_1(f).
\end{equation*}
for any $x\in X_f(\overline{K})$.
\end{theorem}
\begin{proof}
Write $n=\dim X$ and $d=\trdeg_{\mathbb{F}_q}K$.
Take a big and nef polarization
$\mathcal{B}=(B;H,\dots,H)$ of $K$ and a model
$(\mathcal{X},\mathcal{L})$ of $(X,L)$ , where 
$\mathcal{L}$ is an ample line bundle on $\mathcal{X}$. Let $\pi:\mathcal{X}\to B$
be the structure morphism. Since the arithmetic degree $\alpha(f,x)$is independent
of the choice of the ample line bundle $\mathcal{L}$. Then we can choose $\mathcal{L}$ arbitrarily.

Denote by $\Delta_x$ the Zariski closure of the image of
$x:\Spec(\overline{K})\to X\to\mathcal{X}$. Then $\Delta_x$ is a integral
closed subscheme of dimension $d$. Let $\mathcal{I}_{\Delta_x}$
be its ideal sheaf. Then we can replace $\mathcal{L}$ by $\mathcal{L}^{\otimes r}$ for $r$ sufficient large such that
$\mathcal{L}\otimes\mathcal{I}_{\Delta_x}$ is globally generated.

We take global sections
$s_1,\dots,s_n\in
    H^0\left(
    \mathcal{X},\mathcal{L}\otimes\mathcal{I}_{\Delta_x}
    \right)$, 
such that
$\dim\left(\ddiv(s_1)\cap\cdots\cap\ddiv(s_n)\right)=d$ and $\Delta_x$ is an irreducible component of this intersection.

Since
\begin{equation*}
    \widehat{\lambda}_1(f)
    =
    \lim\limits_{m\to\infty}
    \widehat{\deg}_{1,\mathcal{L}}(F^m)^{1/m},
\end{equation*}
Then for every $\epsilon>0$, there exists $C_{\epsilon}>0$, such that
\begin{equation*}
    \widehat{\deg}_{1,\mathcal{L}}(F^m)
    \leq
    C_{\epsilon}
    \left(\widehat{\lambda}_1(f)+\epsilon\right)^m
\end{equation*}
for all $m\geq1$.

Now consider the resolution of indeterminacy of $F^m$:
\begin{center}
    \begin{tikzcd}
        & \Gamma_{F^m}
        \arrow{dl}[swap]{\pi_1}
        \arrow{dr}{\pi_2} & \\
        \mathcal{X}
        \arrow[dashed]{rr}{F^m}
        && \mathcal{X}.
    \end{tikzcd}
\end{center}
On one hand , we know that the strict transformation $\pi_1^{\#}\Delta_x$  of $\Delta_x$ is an irreducible component of $\ddiv(\pi_1^{*}s_1)\cap\cdots\cap\ddiv(\pi_1^{*}s_n)$. On the other hand, consider the height of point $f^m(x)$
\begin{equation*}
    \begin{array}{rl}
      h^{\mathcal{B}}_{\mathcal{X},\mathcal{L}}(f^m(x)) :=\dfrac{\left(\mathcal{L}\cdot \pi^{*}H^{d-1}\cdot \Delta_{f^m(x)}\right)}{[K\left(f^m\left(x\right)\right):K]}=\dfrac{\left(\pi_2^{*}\mathcal{L}\cdot \pi_2^{*}\pi^{*}H^{d-1}\cdot \pi_1^{\#}\Delta_x \right)}{[K(x):K]}
    \end{array}
\end{equation*}
Now we apply Lemma \ref{inequality between line bundles and subvarieties} , one gets
\begin{equation*}
            \left(\pi_2^{*}\mathcal{L}\cdot \pi_2^{*}\pi^{*}H^{d-1}\cdot \pi_1^{\#}\Delta_x \right)\leq\left(\pi_2^{*}\mathcal{L}\cdot \pi_2^{*}\pi^{*}H^{d-1}\cdot\pi_1^{*}\mathcal{L}^n\right)
            =\widehat{\deg}_{1,\mathcal{L}}(F^m)
            \leq C_{\epsilon}(\widehat{\lambda}_1(f)+\epsilon)^m
\end{equation*}
Thus
\begin{equation*}
    \begin{array}{rl}
          \overline{\alpha}_f(x)&=\limsup\limits_{m\to\infty}h^{\mathcal{B}}_{\mathcal{X},\mathcal{L}}(f^m(x))^{1/m}\\
          &=\limsup\limits_{m\to\infty}\dfrac{\left(\pi_2^{*}\mathcal{L}\cdot \pi_2^{*}\pi^{*}H^{d-1}\cdot \pi_1^{\#}\Delta_x \right)^{1/m}}{[K(x):K]^{1/m}}\\
          &\leq \limsup\limits_{m\to\infty}\left( \dfrac{C_{\epsilon}}{[K(x):K]}\right)^{1/m} \left(\widehat{\lambda}_1(f)+\epsilon\right)\\
          &=\widehat{\lambda}_1(f)+\epsilon
    \end{array}
\end{equation*}
Let $\epsilon\to 0$ and applying Theorem \ref{relation between arithmetic dynamical degree and dynamical degree}, we obtain the fundamental inequality
\begin{equation*}
    \overline{\alpha}_f(x)\leq \widehat{\lambda}_1(f)=\max\left\{\lambda_1(f),1\right\}=\lambda_1(f)
\end{equation*}
\end{proof}

Motivated by the results above, we formulate the
Kawaguchi--Silverman conjecture over finitely generated fields
of positive characteristic.
\begin{conjecture}[The Kawaguchi--Silverman conjecture in positive characteristic]\label{KSC}
Let $K$ be a finitely generated field over $\mathbb{F}_q$,
$X$ be a normal projective variety over $K$, 
and $f:X\dashrightarrow X$ be a dominant rational self-map.
Then
\begin{itemize}
  \item the limit $\alpha_f(x)$ in
    Definition \ref{Arithmetic degree} exists for every
    $x\in X_f(\overline{K})$;
    \item for any $x\in X_f(\overline{K})$, if the $f$-orbit
    of $x$ is Zariski dense in $X$, then
    $\alpha_f(x)=\widehat{\lambda}_1(f)$.
\end{itemize}
\end{conjecture}

\begin{theorem}\label{KSC with special Divisor}
 Let $f:X\to X$ be an endomorphism. Suppose there is a big $\mathbb{R}$- Cartier divisor D, such that $f^{*}D\equiv \lambda_1(f) D$. Then the conjecture \ref{KSC} holds.
\end{theorem}
\begin{proof}
Let $H$ be an ample Cartier divisor
on $X$ defined over $K$. Since $X$ is normal and $K$ is
finitely generated over its prime field, $\operatorname{Pic}(X)$
is finitely generated \cite[Proposition~6.1]{GJRW1996}. So there is a monic polynomial $P_f[t]\in \mathbb{Z}[t]$ with $P_f(f^{*})(D)\sim 0$ for all $D\in \Pic(X)$. Let
\begin{equation*}
    V=\Span_{\mathbb{Q}}\{H,f^{*}H,\dots,(f^{*})^{\rho-1}H\}\subset \Pic(X)_{\mathbb{Q}},
\end{equation*}
where $\rho=\deg~P[t]$. Since $P_f(f^*)(H)=0$, the subspace $V$ is $f^*$-invariant. By the same Jordan block height argument in \cite[Theorem~2]{KS2016}, we see that for every
$x\in X(\overline K)$, the arithmetic degree $\alpha_f(x)$ exists and it either equal to 1 or  the absolute value of one of the eigenvalues of $f^{*}$ acting on the finite dimensional vector space $V$.

Now we prove the second part. 
Suppose $\lambda_1(f)=1$, in this case 
  \begin{equation*}
    1\leq \underline{\alpha}_f(x)\leq \overline{\alpha}_f(x)\leq \widehat{\lambda}_1(f)=1
  \end{equation*}
The consequence follows immediately. So we may assume $\widehat{\lambda}_1>1$.

Since $D$ is a big divisor, there exist an ample divisor $A$ and an effective divisor $E$ with $D\sim A+E$. 
 Let $x\in X_f(\overline{K})$ be an point. By fundamental inequality, it suffces to prove $\overline{\alpha}_f(x)\geq \widehat{\lambda}_1(f)$ . Applying proposition \ref{height rela } to $\widehat{\lambda}_1 D-f^{*}D$ and ample divisor $A$, we obtain
\begin{equation*}
    h_{\widehat{\lambda}_1 D-f^{*}D}(x)\leq C\sqrt{h_A(x)}.
\end{equation*}
which menas 
\begin{equation*}
    h_D\left(f(x)\right)\geq \lambda_1h_D(x)-C\sqrt{h_A(x)}.
\end{equation*}
  Iterating this inequality, one gets
  \begin{equation*}
    h_D(f^n(y))
    \geq
    \lambda^n h_D(y)
    -
    C\sum_{i=0}^{n-1}
    \lambda^{n-1-i}\sqrt{h_A(f^i(y))}.
\end{equation*}
One the other hand, since $f$ is a morphismt, he argument of  \cite[Theorem~5.2.1, Remark~5.2.2]{Ohnishi2022} applies
in our setting, which our Proposition \ref{height rela }  bound
the heights associated with numerically trivial divisors. Thus, for every
$\epsilon>0$, there exists a constant $C_{\epsilon}>0$ such that
\begin{equation*}
    h_A(f^n(x))
    \leq C_{\epsilon}(\lambda+\epsilon)^n h_A(x).
\end{equation*}
for all $n\geq0$ and $y\in X(\overline{K})$. Therefore, we can choose  
 $0<\epsilon<\lambda^2-\lambda$, and set
$
    q:=\frac{\sqrt{\lambda+\epsilon}}{\lambda}<1.
$
Then
\begin{equation*}
    \begin{aligned}
        h_D(f^n(x))
        &\geq
        \lambda^n
        \left(
        h_D(x)
        -
        \frac{C_0\sqrt{C_{\epsilon}}}{\lambda}
        \sum_{i=0}^{n-1}q^i\sqrt{h_A(x)}
        \right)\\
        &\geq
        \lambda^n
        \left(h_D(x)-C_1\sqrt{h_A(x)}\right),
    \end{aligned}
\end{equation*}
where $C_1:=\frac{C_0\sqrt{C_{\epsilon}}}{\lambda(1-q)}$ is independent of $n$ and $x$.

Since $\mathcal{O}_f(x)$ is Zariski dense, by Northcott property, the set $\{h_{A}(f^n(x))\}_{n\geq 0}$ is not bounded from above. Thus, we may assume $h_{A}(x)> C_1^2$ and $x\in \left(X\setminus \supp(E)\right)(\overline{K})$. Then
\begin{equation*}
  h_D(x)-C_1 \sqrt{h_A(x)}\geq h_A(x)-C_1\sqrt{h_A(x)}>0.
\end{equation*}
it shows that there exists a $c>0$ such that
\begin{equation*}
    h_D(f^n(x))\geq c\widehat{\lambda}_1(f)^n.
\end{equation*}
taking $n$-th root of both sides, and letting $n\to \infty$, we get
  \begin{equation*}
   \widehat{\lambda}_1(f) \geq\overline{\alpha}_f(x)\geq\underline{\alpha}_f(x)\geq\widehat{\lambda}_1(f).
\end{equation*}
which concludes the proof.
\end{proof}

The same estimates also give the following corollary.

\begin{corollary}\label{canonical height}
Let $f:X\to X$ be a dominant endomorphism defined over $K$,
and let $D$ be a nef $\mathbb{R}$-Cartier divisor such that
\begin{equation*}
    f^*D\equiv\lambda_1(f)D,
    \qquad \lambda_1(f)>1.
\end{equation*}
Then, for every $x\in X(\overline{K})$, the limit
\begin{equation*}
    \widehat{h}_D(x)
    :=
    \lim\limits_{n\to\infty}
    \frac{h_D(f^n(x))}{\lambda_1(f)^n}
\end{equation*}
exists. We call $\widehat{h}_D$ the canonical height
with respect to $D$. Moreover,
\begin{equation*}
    \widehat{h}_D(x)
    =h_D(x)+O\left(\sqrt{h_A(x)}\right),
\end{equation*}
where $A$ is an ample divisor and $h_A\geq1$.
\end{corollary}

\section{ Arithmetic degrees over fields of positive characteristic}\label{algebraic closed case}
In this section, we study arithmetic degrees over a algebraically closed fields of positive characteristic. Our approach rests on
the fact that a dynamical system, together with a chosen
point, can be defined over a finitely generated subfield.
This allows us to apply the height theory developed above
to study the growth of heights along its orbit.
\subsection{Fields of definition and arithmetic degrees}

Let $\mathbf{k}$ be an algebraically
closed field of characteristic $p>0$. Let X be a normal projective variety over $\mathbf{k}$
 and $f:X\dashrightarrow X$ be a dominant rational
self-map over $\mathbf{k}$.
We extend the notion of arithmetic degree to this setting
by working over finitely generated fields of definition.

 More precisely, For any $x\in X_f(\mathbf{k})$, the triple $(X,f,x)$
is determined by finitely many equations and coordinates.
Thus we can take a finitely
generated field
$K/k$. such that there is  a normal
projective variety $X_K$ over $K$, a dominant rational self-map
$f_K:X_K\dashrightarrow X_K$, and a point $x_K\in X_K(K)$,
whose base change to $\mathbf{k}$ recovers $(X,f,x)$, i.e., there is an isomorphism
\begin{equation*}
    X\simeq X_K\times_{\operatorname{Spec}K}
    \operatorname{Spec}\mathbf{k}.
\end{equation*}
which $f$ and $x$ are identified with the base changes
of $f_K$ and $x_K$, respectively.
We call the triple $(X,f,x)$ can be defined over $K$.
\[
\begin{tikzcd}
    & X \arrow{dd} \arrow[dashed]{rr}{f} \arrow{rd} & & X \arrow{dd} \arrow{ld} \\
    \Spec\mathbf{k} \arrow{dd} \arrow{rr}{\quad\quad\mathrm{id}} \arrow{ru}{x} & & \Spec\mathbf{k} \arrow{dd} & \\
    & X_K \arrow{rd} \arrow[dashed]{rr}{\quad \quad f_K} & & X_K \arrow{ld} \\
    \Spec K \arrow{ru}{x_K} \arrow{rr}{\mathrm{id}} & & \Spec K &
\end{tikzcd}
\]

The next proposition proves that the 
arithmetic degrees are independent of this choice.
This allows us to study dynamical system
$(X,f,x)$ over $\mathbf{k}$.

\begin{proposition}
Let $(X,f,x)$ be as above.
For any finitely generated field $K\subseteq\mathbf{k}$ 
over $k$ which $(X,f,x)$ can be defined,
denote the corresponding system by $(X_K,f_K,x_K)$.
Then the arithmetic degrees
$\overline{\alpha}_{f_K}(x_K)$ and
$\underline{\alpha}_{f_K}(x_K)$
are independent of the choice of $K$.
\end{proposition}

\begin{proof}
Suppose  $(X,f,x)$ can be defined over two finitely
generated fields $K_1,K_2\subseteq\mathbf{k}$,
with corresponding triples $(X_{K_i},f_{K_i},x_{K_i})$.
If the orbit of $x$ is finite, both arithmetic degrees
are equal to $1$, and there is nothing to prove.
So We may assume that the orbit is infinite.

Again choose a finitely generated field $K_3$, 
containing $K_1$ and $K_2$
we choose a polarization $\mathcal{B}_i$
of $K_i$ and a model $(\mathcal{X}_i,\mathcal{L}_i)$
of $X_{K_i}$  with $\mathcal{L}_i$ ample. Denote
$  h_i
    :=h^{\mathcal{B}_i}_{(\mathcal{X}_i,\mathcal{L}_i)}.$ 
Applying Theorem~\ref{comparable of height under field extension}
to $K_3/K_1$ and $K_3/K_2$, then there exists constants
$N_1,N_2>0$ and $M_1,M_2\in\mathbb{R}$ such that
\begin{equation*}
N_1h_1\left(f_{K_1}^{n}(x_{K_1})\right)+M_1
    \leq
    h_2\left(f_{K_2}^{n}(x_{K_2})\right)
    \leq
    N_2h_1\left(f_{K_1}^{n}(x_{K_1})\right)+M_2
\end{equation*}
for every $n\geq0$, which the constants are independent of $n$.
Taking $n$-th roots and then the upper and lower limits, we obtain
\begin{equation*}
    \overline{\alpha}_{f_{K_1}}(x_{K_1})
    =
    \overline{\alpha}_{f_{K_2}}(x_{K_2}),
    \qquad
    \underline{\alpha}_{f_{K_1}}(x_{K_1})
    =
    \underline{\alpha}_{f_{K_2}}(x_{K_2}).
\end{equation*}
which means the arithmetic degrees are well-defined.
\end{proof}

By this proposition, we can define the  arithmetic degrees of any $x\in X_f(\mathbf{k})$
by doing this on a finitely generated field.
\begin{definition}[Arithmetic degree]
Let $X$ be a normal projective variety over
$\mathbf{k}$ and 
$f:X\dashrightarrow X$ be a dominant rational self-map.
For any $x\in X_f(\mathbf{k})$, we can choose a finitely generated
$K$ over $\mathbb{F}_q$, such that  $(X,f,x)$ can be defined. Then, we define the  arithmetic degrees of $x$ by
\begin{equation*}
    \overline{\alpha}_f(x)
    :=\overline{\alpha}_{f_K}(x_K),
    \qquad
    \underline{\alpha}_f(x)
    :=\underline{\alpha}_{f_K}(x_K).
\end{equation*}
\end{definition}

\begin{theorem}\label{fundemental inequality}
For every $x\in X_f(\mathbf{k})$, we have
\begin{equation*}
    \overline{\alpha}_f(x)\leq\lambda_1(f).
\end{equation*}
\end{theorem}

\begin{proof}
Choose a finitely generated field
$K\subseteq\mathbf{k}$ such that $(X,f,x)$ is defined,
and denote the corresponding triple by $(X_K,f_K,x_K)$.
Since dynamical degrees are invariant under
base change, we have
\begin{equation*}
    \lambda_1(f_K)=\lambda_1(f).
\end{equation*}
If $K$ is finite, the orbit of $x_K$ is finite,
so $\overline{\alpha}_f(x)=1\leq\lambda_1(f)$.
Otherwise, by Theorem
\ref{fundemental inequality over fg case}
we obtain
\begin{equation*}
    \overline{\alpha}_f(x)
    =
    \overline{\alpha}_{f_K}(x_K)
    \leq
    \widehat{\lambda}_1(f_K)
    =
    \lambda_1(f_K)
    =
    \lambda_1(f).
\end{equation*}
\end{proof}

These results lead us to formulate the
Kawaguchi--Silverman conjecture in positive characteristic setting.

\begin{conjecture}
Let $X$ be a normal projective variety over a field
$\mathbf{k}$ 
and  $f:X\dashrightarrow X$ be a dominant
rational self-map. Then 
\begin{itemize}
    \item
    The arithmetic degree $\alpha_f(x)$ exists
    for every $x\in X_f(\mathbf{k})$.

    \item
    For any $x\in X_f(\mathbf{k})$, if the $f$-orbit $O_f(x)$ is Zariski dense,
   then $
        \alpha_f(x)=\lambda_1(f).$
\end{itemize}
\end{conjecture}

By spreading out to a finitely generated subfield,
we can extend several results established in the previous section
to $\mathbf{k}$. 
In particular, we obtain the following cases in which
the Kawaguchi--Silverman conjecture holds.
\begin{proposition}
Let $X$ be a  normal  projective variety over $\mathbf{k}$
and $f:X\to X$ be a endomorphism.
Suppose that there is a big $\mathbb{R}$-Cartier
divisor $D$ on $X$ such that
\begin{equation*}
    f^*D\equiv\lambda_1(f)D.
\end{equation*}
Then the Kawaguchi--Silverman conjecture holds for $f$.
\end{proposition}

\begin{proof}

Choose a finitely generated subfield
$K\subseteq\mathbf{k}$, such that the system $(X,f,x)$
and the  divisor $D$ can be defined.
Denote the corresponding divisor by $D_K$. Then  $D_K$ is big and
\begin{equation*}
    f_K^*D_K\equiv\lambda_1(f_K)D_K.
\end{equation*}
Applying Theorem \ref{KSC with special Divisor}, we obtain the arithmetic degree $\alpha_{f_K}(x_K)$ exists. Moreover, if $\mathcal{O}_f(x)$ is Zariski dense, one gets
\begin{equation*}
    \alpha_f(x)
    =
    \alpha_{f_K}(x_K)
    =
    \lambda_1(f_K)
    =
    \lambda_1(f).
\end{equation*}
This completes the proof.
\end{proof}
We next study the existence and possible values of arithmetic
degrees for endomorphisms.
Over fields of characteristic zero, Luo and Song
\cite{LS2025} proved that these degrees exist and take
only finitely many values, all of which are algebraic integers,
by reducing to finitely generated fields over $\mathbb{Q}$.

In positive characteristic, we instead work over fields
finitely generated over $\mathbb{F}_p$ and use the
heights constructed in the Section~\ref{height chatper}.
The associated proper adelic structures allow us
to apply \cite[Theorem~5.3.1]{Ohnishi2022},
which imposes no restriction on the characteristic.

\begin{proposition}
Let $X$ be a normal projective variety over $\mathbf{k}$
and $f:X\to X$ be an endomorphism.
Then
\begin{itemize}
    \item
    The arithmetic degree $\alpha_f(x)$ exists
    for every $x\in X(\mathbf{k})$.

    \item
    The arithmetic degree $\alpha_f(x)$
    is an algebraic integer.

    \item
    The set
    \begin{equation*}
        \left\{
            \alpha_f(x)\mid x\in X(\mathbf{k})
        \right\}
        \subseteq\mathbb{R}
    \end{equation*}
    is finite.
\end{itemize}
\end{proposition}

\begin{proof}

  For any $x \in X(\mathbf{k})$. We can choose a finitely generated field $K$ such that $(X,f,x)$ can be defined. By definition, $\alpha_f(x)=\alpha_{f_{K}}(x_K)$. Choose an ample polarization $\mathcal{B}$ of $K$,
and let $S_{K,\mathcal{B}}$ be the associated proper adelic curve.  Now we apply \cite[Theorem~5.3.1]{Ohnishi2022} to $(X_K,f_K)$ with respect to this adelic structure, one gets
 \begin{itemize}
    \item
    The arithmetic degree $\alpha_{f_K}(x_K)$ exists
    for every $x_K\in X(K)$.

    \item
    The arithmetic degree $\alpha_{f_K}(x_K)$ 
    is an algebraic integer.

    \item
    The set
    \begin{equation*}
        \left\{
            \alpha_{f_K}(x_K) \mid x_K \in X(K)
        \right\}
        \subseteq\mathbb{R}
    \end{equation*}
    is finite.
\end{itemize}

  We now using these  to prove the finiteness of
$ \left\{
            \alpha_f(x)\mid x\in X(\mathbf{k})
        \right\}$.
We want to control the arithmetic degree by considering the
action of $f^*$ on $\NS(X)_{\mathbb{Q}}$.

  We know that the Nron--Severi group $\NS(X)_{\mathbb{Q}}$ is finitely generated. Choose $D_1,\dots,D_k$ be a basis of $\NS(X)_{\mathbb{Q}}$, and a finitely generated field $K^{'}\supseteq K$ such that all divisors $D_i$ and $x$ can be defined. Under these assumptions, we have a canonical isomorphism $\NS(X_{K^{'}})_{\mathbb{Q}}\xrightarrow{\sim}\NS(X)_{\mathbb{Q}}$ .
Again by the arguments in \cite[Theorem~5.3.1]{Ohnishi2022}, we know that the arithmetic degree $\alpha_{f_{K^{'}}}(x_{K^{'}})$ equal to 1 or the absolute values of  some eigenvalue of $f^{*}_{K^{\prime}}$, which is finite. Moreover, the isomorphism  commutes with the action of $f^{*}_{K^{'}}$ and $f^{*}$. Consider the following diagram:
  \begin{center}
    \begin{tikzcd}
        \NS(X_{K'})_{\mathbb{Q}} \arrow{r}{\sim} \arrow{d}{f_{K'}^*} & \NS(X)_{\mathbb{Q}} \arrow{d}{f^*} \\
        \NS(X_{K'})_{\mathbb{Q}} \arrow{r}{\sim} & \NS(X)_{\mathbb{Q}}
    \end{tikzcd}
\end{center}
So the eigenvalues of $f^{*}$ is exactly the eigenvalues of $f^{*}_{K'}$. Thus, the set
\begin{equation*}
    \{\alpha_f(x)~|~x\in X(\mathbf{k})\}\subset \{1\}\cup \{|\lambda|~|~\lambda \text{ is an eigenvalue of }
f^{*}:\NS(X_{K^{\prime}})_{\mathbb{Q}}\to \NS(X)_{\mathbb{Q}}\}
\end{equation*}
is finite, so we completes the proof. 
\end{proof}

\subsection{Arithmetic degrees of transcendental points}
Throughout this subsection, we assume that \(\mathbf{k}\) is an algebraically closed field of characteristic $p$ with $\trdeg_{\mathbb{F}_p}(\mathbf{k})=\infty$.  We follow the terminology of \cite[Definition~0.3]{LS2025}.

\begin{definition}
  Let  $X_K$ be a variety over a subfield $K\subset \mathbf{k}$ and $X=X_K\times_{K}\Spec \mathbf{k}$.  Denote $\pi_K: X\to X_K$ the canonical projection. For a point $x\in X(\mathbf{k})$ , denote $Z(x)_K$ the Zariski closure of $\pi_K(x)$ in $X_K$ and $Z(x)=\pi_K^{-1}(Z(x)_K)$. We call $\dim Z(x)$ the transcendence degree of x with respect to $X_K\to \Spec K$. We say x is transcendental if $\dim Z(x)\geq 1$ and call it purely transcendental if $\dim Z(x)=\dim X$.
\end{definition}

So a point $x$ is purely transcendental if and only if it does not lie in any proper closed $K$-subvarieties of X. Denote
\begin{equation*}
  S:=\{V\subsetneq X \mid V\text{is a closed subvariety defined over K}\}
\end{equation*}
Thus, the set of purely transcendental points with respect to $X_K\to \Spec K$ is given by $X(\mathbf{k})\setminus \bigcup\limits_{V\in S} V(\mathbf{k})$.

\begin{proposition}\label{galois action on purely transcendental point}
    Let $K\subseteq\mathbf{k}$ be a finitely generated field over
$\mathbb{F}_p$, and let $X_K$ be any variety over $K$. Put $X=X_K\times_K\Spec\mathbf{k}$.
 Then for any purely transcendental points $x,y\in X(\mathbf{k})$ with respect to $X_K\to \Spec K$, there exists a $K$-automorphism $\sigma\in\Gal(\mathbf{k}/K)$ with $y=\sigma(x)$.
\end{proposition}
\begin{proof}
In \cite{LS2025}, Luo and Song  proved the 
result for quasi-projective varieties over fields finitely generated
over $\mathbb{Q}$. We give an alternative proof
that also for arbitrary varieties.

Let $\pi_K:X\to X_K$ be the canonical projection.
Since $x$ and $y$ are purely transcendental with
respect to $X_K/K$, both $\pi_K(x)$ and $\pi_K(y)$
are the generic point of $X_K$.
Choose a nonempty affine open subset
$U_K=\Spec(R)\subseteq X_K$.
Then $x,y\in U_K(\mathbf{k})$, so we get two
 $K$-embeddings
\begin{equation*}
    \iota_x,\iota_y:K(X_K)=\kappa(\eta)\hookrightarrow\mathbf{k}.
\end{equation*}
Denote their images by $K_x$ and $K_y$, respectively. Then
$\tau=\iota_y\circ\iota_x^{-1}:K_x\to K_y$ is a $K$-isomorphism.

Choose transcendence bases $T_x$ and $T_y$ of $\mathbf{k}$ over
$K_x$ and $K_y$, respectively. Since $K_x$ and $K_y$ are finitely
generated over $\mathbb{F}_p$ and $\trdeg_{\mathbb{F}_p}\mathbf{k}$
is infinite, the sets $T_x$ and $T_y$ have the same cardinality. Choose a bijection
$\varphi:T_x\to T_y$. Applying $\tau$ to the coefficients and
sending each $t\in T_x$ to $\varphi(t)$ gives an  isomorphism $ \widetilde{\tau}:K_x(T_x)\xrightarrow{\sim}K_y(T_y).$

Moreover, 
since $\mathbf{k}$ is an algebraic closure of both fields,
$\widetilde{\tau}$ extends to a $K$-automorphism $\sigma$ of
$\mathbf{k}$ by \cite[Tag~09GV]{stacks-project}.
For every $a\in R$, we have
\begin{equation*}
    \sigma(a(x))=\sigma(\iota_x(a))=\iota_y(a)=a(y).
\end{equation*}
Thus $\sigma(x)=y$, which completes the proof.
\end{proof}

\begin{theorem}\label{Arithmetic degree under galois action}
Let $X$ be a normal quasi-projective variety over
$\mathbf{k}$, and let $f:X\dashrightarrow X$ be a
dominant rational self-map.
Suppose that $(X,f)$ is defined over a subfield
$K\subseteq\mathbf{k}$ finitely generated over
$\mathbb{F}_p$.
For any $x,y\in X_f(\mathbf{k})$, if there exists
$\sigma\in\Gal(\mathbf{k}/K)$ such that $\sigma(x)=y$,
then
\begin{equation*}
    \overline{\alpha}_f(x)=\overline{\alpha}_f(y),
    \qquad
    \underline{\alpha}_f(x)=\underline{\alpha}_f(y).
\end{equation*}
\end{theorem}

\begin{proof}
Choose a finitely generated extension
$K_1/K$ contained in $\mathbf{k}$ over which
$(X,f,x)$ is defined.
Then $(X,f,y)$ is defined over $K_2=\sigma(K_1)$,
and $\sigma:K_1\to K_2$ is a $K$-isomorphism.

Fix an ample polarization
$\mathcal{B}_1=(B_1;H_1,\dots,H_{d-1})$ of $K_1$,
where $d=\trdeg_{\mathbb{F}_p}K_1$.
Via the isomorphism, it's also an ample polarization of $K_2$.
Take an embedding $\mathcal{X}\hookrightarrow \mathbb{P}^N$, 
For $P=[a_0:\cdots:a_N]\in X_K(K_1)$,
consider the naive height
\begin{equation*}
\begin{aligned}
    h_{\mathrm{nv},K_2}^{\mathcal{B}_2}(\sigma(P))
    &=
    \sum_{\Gamma\in B_1^{(1)}}
    \max_i\{-\ord_{\Gamma}(\sigma(a_i))\}
    \left(\Gamma\cdot
    H_1\cdots H_{d-1}\right)\\
    &=
    \sum_{\Gamma\in B_1^{(1)}}
    \max_i\{-\ord_\Gamma(a_i)\}
    \left(\Gamma\cdot H_1\cdots H_{d-1}\right)\\
    &=
    h_{\mathrm{nv},K_1}^{\mathcal{B}_1}(P).
\end{aligned}
\end{equation*}
Similarly, we have
$    h_{\mathrm{nv},K_1}^{\mathcal{B}_1}(f^n(x))
    =
    h_{\mathrm{nv},K_2}^{\mathcal{B}_2}(f^n(y))$ for every $n\geq 0$.
taking $n$-th roots, and then taking the upper and
lower limits, then we obtain
$\overline{\alpha}_f(x)=\overline{\alpha}_f(y)
    $ and $
    \underline{\alpha}_f(x)=\underline{\alpha}_f(y)$.

\end{proof}

\begin{corollary}
  Let $X$ be a normal projective variety and $f:X\dashrightarrow X$ be a dominant rational self-map over $\mathbf{k}$. Suppose $(X,f)$ is defined over a finitely generated field $K$. For any $x,y\in X_f(\mathbf{k})$ satisfying $Z(x)_{K}=Z(y)_{K}$, we have
  \begin{equation*}
    \overline{\alpha}_f(x)=\overline{\alpha}_f(y)~~ and~~
    \underline{\alpha}_f(x)=\underline{\alpha}_f(y)
  \end{equation*}
In particular, there exists a countable set of proper subvarieties of $X$:
\begin{equation*}
  S=\{V\subsetneq X~|~V ~\text{is a closed subvariety defined over K}\}\cup\{I(f^n)~|~n\in \mathbb{Z}_{+}\}
\end{equation*}
And for any $x,y\in X(\mathbf{k})\setminus \bigcup\limits_{V\in S}V(\mathbf{k})$, we have $ \overline{\alpha}_f(x)=\overline{\alpha}_f(y)$ and  
   $ \underline{\alpha}_f(x)=\underline{\alpha}_f(y)$.
\end{corollary}
\begin{proof}
  Since $(X,f)$ is defined over $K$, there is a normal
projective variety $X_K$ over $K$ such that $
    X=X_K\times_{\Spec(K)}\Spec(\mathbf{k})$.
Set $Y=Z(x)$ and $Y_K=Z(x)_K$.
Then $x$ and $y$ are both purely transcendental points with respect to $Y_K\to\Spec(K)$.
By Proposition \ref{galois action on purely transcendental point},
there exists a $K$-automorphism
$\sigma\in\Gal(\mathbf{k}/K)$ such that $\sigma(x)=y$. Applying Theorem \ref{Arithmetic degree under galois action},
we obtain the desired equalities.

Furthermore, let
\begin{equation*}
    S_1=
    \{V\subsetneq X\mid
    V\text{ is a closed subvariety defined over }K\},
    \qquad
    S_2=\{I(f^n)\mid n\in\mathbb{Z}_{+}\}.
\end{equation*}
Since $K$ is countable, $S_1$ is countable.
Thus $S=S_1\cup S_2$ is also countable.
On the other hand, we know $
    X(\mathbf{k})\smallsetminus
    \bigcup_{V\in S}V(\mathbf{k})$ 
equal to the set of purely transcendental points
with respect to $X_K\to\Spec(K)$.
For any two such points $x$ and $y$, we have
$Z(x)_K=Z(y)_K=X_K$.
Hence their  arithmetic degrees
are equal.
\end{proof}

For a fixed finitely generated field $K$, there are at most
countably many closed subvarieties $Z(x)_K$ of $X_K$.
We therefore immediately obtain the following corollary.

\begin{corollary}
With the notation above, the set
\begin{equation}
    \{\overline{\alpha}_f(x)\mid x\in X_f(\mathbf{k})\}
    \cup
    \{\underline{\alpha}_f(x)\mid x\in X_f(\mathbf{k})\}
\end{equation}
is at most countable.
Furthermore, if $f$ is an endomorphism, then this set
is a finite set of algebraic integers.
\end{corollary}

\begin{lemma}\label{finiteness of transcendental points}
Let $X$ be a quasi-projective variety
over $\mathbf{k}$, defined over a finitely generated subfield
$K\subseteq\mathbf{k}$.
For any $\beta\in\mathbf{k}\smallsetminus\overline{K}$,
put $L=K(\beta)$.
Then there exists a positive integer $N$ such that the set
  \begin{equation*}
        \{x\in X_L(\overline{L})~|~\text{$x$ is transcendental with respect to $X_K\to \Spec(K)$, $[L(x):L]<N$}\}
    \end{equation*}
is infinite.
\end{lemma}

\begin{proof}
Replacing $X_K$ by a nonempty affine open subset,
we may assume that
\begin{equation*}
    X_K=\Spec(R_K),
    \qquad
    R_K=K[x_1,\dots,x_n]/(f_1,\dots,f_m).
\end{equation*}
Write $d=\dim X\geq1$.
By Noether normalization, there exists a finite surjective
morphism
\begin{equation*}
    f_K:X_K\longrightarrow\mathbb{A}^d_K,
\end{equation*}
which induces a finite surjective morphism
$f:X\to\mathbb{A}^d$ by base change.
There exists an integer $N>0$ such that the set
  \begin{equation*}
        \{x\in \mathbb{A}^d(L)~|~\text{$x$ is transcendental with respect to $\mathbb{A}^d_K\to \Spec(K)$}\}
    \end{equation*}
is infinite. On the other hand, since the morphism $f_K$ is finite, $R_L=L[x_1,\dots,x_n]/(f_1,\dots,f_m)$ is finitely generated $L[y_1,\dots,y_d]$-module. Suppose it can be generated by $r$ elements. Then for any $x\in X_L(\overline{L})$, with $y=f_L(x)$. The coordinate ring of the fiber over $y$
has dimension at most $r$ over $L(y)$,
and $L(x)$ is a residue field of this fiber. Then, we have
\begin{equation*}
    [L(x):L]\leq r[L(y):L].
\end{equation*}
 Furthermore, if y is transcendental with respect to $\mathbb{A}_K^d\to \Spec K$, it follows that $x$ is transcendental with respect to $X_K\to \Spec K$. Thus, the set
 \begin{equation*}
  \{x\in X_L(\overline{L})~|~x~ \text{is transcendental with respect to $X_K\to \Spec K$, $[L(x):L]\leq r$}\}
 \end{equation*}
 is infinite.
\end{proof}

\begin{theorem}
Let $X$ be a normal projective variety over $\mathbf{k}$ and let $f:X\to X$ be a surjective
endomorphism.
Suppose that $(X,f)$ can be defined over a finitely generated field $K\subseteq\mathbf{k}$.
Then, for any purely transcendental point
$x\in X(\mathbf{k})$ with respect to
$X_{\overline{K}}\to\Spec(\overline{K})$, we have
$\alpha_f(x)=\lambda_1(f)$.
In particular, there exists a countable collection $S$
of proper closed subvarieties of $X$ such that
$\alpha_f(x)=\lambda_1(f)$
for every
$x\in X(\mathbf{k})\smallsetminus
\bigcup_{V\in S}V(\mathbf{k})$.
\end{theorem}
\begin{proof}
  We may assume that $\lambda_1(f)>1$, since the case
$\lambda_1(f)=1$ follows from the fundamental inequality.

Suppose $(X,f)$ can be defined over a finitely generated field $K$. Let 
\begin{equation*}
    S:=
    \left\{
        V\subsetneq X
        \mid V\text{ is a closed subvariety of $X$ defined over }K
    \right\}.
\end{equation*}
Since $K$ is countable, the collection $S$ is countable.
The complement
\begin{equation*}
    X(\mathbf{k})\smallsetminus
    \bigcup_{V\in S}V(\mathbf{k})
\end{equation*}
is the purely transcendental points
with respect to $X_K\to\Spec(K)$.
Since $\overline{K}/K$ is algebraic, these are also precisely
the purely transcendental points with respect to
$X_{\overline{K}}\to\Spec(\overline{K})$.
By the preceding result, we have
\begin{equation*}
    \alpha_f(x)=\alpha_f(y)
\end{equation*}
for any two points $x,y$ in this complement.
Thus it suffices to find one point
\begin{equation*}
    x_0\in X(\mathbf{k})\smallsetminus
    \bigcup_{V\in S}V(\mathbf{k})
\end{equation*}
such that $\alpha_f(x_0)=\lambda_1(f)$.

First,choose ample effective Cartier divisors
$H_1,\dots,H_{d-1}$ on $X$, Applying s Bertini theorem
\cite[Theorem~$7'$]{Seidenberg1950} d-1 times, we get get an we obtain a normal
\begin{equation*}
    C:=H_1\cap\cdots\cap H_{d-1}
\end{equation*}
and it is not contained in any $V\in S$.

Since $f^*$ preserves the nef cone and has spectral radius
$\lambda_1(f)$ on $\NS(X)_{\mathbb{R}}$,
there exists a nef $\mathbb{R}$-divisor D such that $f^{*}D\equiv \lambda_1(f)D$. Then we have $\left(D\cdot H_1\cdots H_{d-1}\right)>0$ by \cite[Lemma 2.5]{Hu2020}. Thus, $D|_{C}$ is ample.

Choose a finitely generated subfield
$K_1\subseteq\mathbf{k}$ containing $K$
over which $X,f,C,D$ are defined.
For every $V\in S$, the intersection $C\cap V$
is zero-dimensional and defined over $K_1$.
Such that $ C(\mathbf{k})\cap\bigcup_{V\in S}V(\mathbf{k})$. Suppose 
\begin{equation*}
  C(\mathbf{k})\cap\bigcup_{V\in S}V(\mathbf{k}) =\{x_1,x_2,\dots\}
\end{equation*}
where $x_1,x_2,\dots$ are defined over $K_2$, which is countably generated extension over $K_1$.  Since $K_2/K_1$ is algebraic and $\trdeg_{K_2}\mathbf{k}=\infty$, we can choose an element $\beta\in \mathbf{k}\setminus\overline{K}_2$ transcendental over $K_2$. Then we have $\overline{K_1(\beta)}\cap \overline{K_2}=\overline{K_1}$. Set $K_3=K_1(\beta)$.

Fix an ample polarization $\mathcal{B}$ of $K_3$
and a corresponding Moriwaki height $h_D$.
For simplicity, we continue to denote $D_{K_3}$ by $D$.
Since $D|_C$ is ample, we may assume that $h_D|_C\geq1$.
Consider the canonical height
\begin{equation*}
    \widehat{h}_D
    :=
    \lim_{n\to\infty}
    \frac{h_D\circ f_{K_3}^n}{\lambda_1(f_{K_3})^n}.
\end{equation*}
By Corollary \ref{canonical height}, the limit exists and
\begin{equation*}
    \widehat{h}_D=h_D+O(\sqrt{h_A}),
\end{equation*}
where $A$ is an ample Cartier divisor and $h_A\geq1$.
Comparing heights on $C$, we obtain constants
$0<C_1<C_2$ such that
\begin{equation*}
    \widehat{h}_D(x)
    \geq h_D(x)-C_1\sqrt{h_A(x)}
    \geq h_D(x)-C_2\sqrt{h_D(x)}
\end{equation*}
for every $x\in C_{K_3}(\overline{K_3})$.

If $\widehat{h}_D(x)>0$, then $\alpha_f(x)=\lambda_1(f)$.
Consequently,
\begin{equation*}
    \begin{array}{rl}
            \{x\in C_{K_3}(\overline{K_3})~|~\alpha_f(x)< \lambda_1(f)\}&\subset \{x\in C_{K_3}(\overline{K_3})~|~\widehat{h}_D(x)\leq 0\}\\
            &\subset \{x\in C_{K_3}(\overline{K_3})~|~h_D(x)-C_2\sqrt{h_D(x)}\leq 0\}\\
            &\subset \{x\in C_{K_3}(\overline{K_3})~|~h_D(x)\leq C_2^2\}
    \end{array}
\end{equation*}
By Lemma \ref{finiteness of transcendental points},
there exists a positive integer $N$ such that the set
 \begin{equation*}
 T= \{x\in C(\overline{K_3})~|~x~ \text{is transcendental with respect to $C_{K_1}\to\Spec K_1$, $[K_3(x):K_3]<N$}\}
 \end{equation*}
is infinite.
Since $D|_C$ is ample, the Northcott property gives
a point $x_0\in T$ such that $h_D(x_0)>C_2^2$.
Thus $\widehat{h}_D(x_0)>0$, and hence
$\alpha_f(x_0)=\lambda_1(f)$.
Moreover, $x_0$ does not belong to
$C(\mathbf{k})\cap\bigcup_{V\in S}V(\mathbf{k})$.
Otherwise, it would be defined over
\begin{equation*}
    \overline{K_2}\cap\overline{K_3}
    =\overline{K_1},
\end{equation*}
contrary to its transcendence over $K_1$.
By the constancy of arithmetic degrees on the complement,
we conclude that
\begin{equation*}
    \alpha_f(x)=\lambda_1(f)
    \qquad
    \text{for every }
    x\in X(\mathbf{k})\smallsetminus
    \bigcup_{V\in S}V(\mathbf{k}).
\end{equation*}

\end{proof}
\section*{AI statement}
This paper is based on the author's master's thesis. While drafting the manuscript, the author used
 GPT and DeepSeek to  help search for known mathematical results and polish the English. After completing
the first draft, the author also used
these tools to check the grammar and
identify possible mathematical errors. The author takes responsibility for the final text and
its mathematical content.

\bigskip
\noindent
\begin{minipage}{\linewidth}
  \small\raggedright
  \hypersetup{urlcolor=black}

  \rule{0.2\linewidth}{0.4pt}\par
  \smallskip
  \textit{\today}\par
  \medskip

  \textsc{Yubo Zhou}, School of Mathematics,
  Harbin Institute of Technology,\\
  \hspace*{2em}150001 Harbin, P.~R.~China
  \quad\textbullet\quad
  \textit{E-mail:}~%
  \href{mailto:1201200209@stu.hit.edu.cn}{\nolinkurl{1201200209@stu.hit.edu.cn}}
\end{minipage}

\begin{thebibliography}{[AB]}

\bibitem[CS93]{CS1993}
G.~S. Call and J.~H. Silverman, Canonical heights on varieties with morphisms, Compositio Math. {\bf 89} (1993), no.~2, 163--205; MR1255693

\bibitem[Cutk15]{Cutk15}
S.~D. Cutkosky, Teissier's problem on inequalities of nef divisors, J. Algebra Appl. {\bf 14} (2015), no.~9, 1540002, 37 pp.; MR3368254

\bibitem[Dang20]{Dang20}
 N.-B. Dang, Degrees of iterates of rational maps on normal projective varieties, Proc. Lond. Math. Soc. (3) {\bf 121} (2020), no.~5, 1268--1310; MR4133708

\bibitem[DS05]{DS05}
T.-C. Dinh and N. Sibony, Une borne sup\'erieure pour l'entropie topologique d'une application rationnelle, Ann. of Math. (2) {\bf 161} (2005), no.~3, 1637--1644; MR2180409

\bibitem[FAG05]{FAG2005}
B. Fantechi, L. G\"ottsche, L. Illusie, S. L. Kleiman, {\it Fundamental algebraic geometry}, Mathematical Surveys and Monographs, 123, American Mathematical Society, Providence, RI, 2005; MR2222646

\bibitem[FOV99]{FOV99}
H. Flenner, L. O'Carroll and W. Vogel, {\it Joins and intersections}, Springer Monographs in Mathematics, Springer, Berlin, 1999; MR1724388

\bibitem[Ful98]{Ful98}
W. Fulton, {\it Intersection theory}, second edition, 
Ergebnisse der Mathematik und ihrer Grenzgebiete. 3. Folge. A Series of Modern Surveys in Mathematics, 2, Springer, Berlin, 1998; MR1644323

\bibitem[EGAIV]{EGA4}
A. Grothendieck, \'El\'ements de g\'eom\'etrie alg\'ebrique. IV. \'Etude locale des sch\'emas et des morphismes de sch\'emas IV, Inst. Hautes \'Etudes Sci. Publ. Math. No. 32 (1967), 361 pp.; MR0238860

\bibitem[GJRW96]{GJRW1996}
R.~M.~Guralnick, D.~B.~Jaffe, W.~Raskind, and R.~Wiegand, On the Picard group: torsion and the kernel induced by a faithfully flat map, J. Algebra {\bf 183} (1996), no.~2, 420--455; MR1399035

\bibitem[Hu20]{Hu2020}F. Hu, A theorem of Tits type for automorphism groups of projective varieties in arbitrary characteristic, Math. Ann. {\bf 377} (2020), no.~3-4, 1573--1602; MR4126902

\bibitem[JSXZ24]{JSXZ2024}
J. Jia, T. Shibata, J. Xie and D.-Q. Zhang,
Endomorphisms of quasi-projective varieties:
towards Zariski dense orbit and Kawaguchi--Silverman conjectures,
Math. Res. Lett. {\bf 31} (2024), no.~3, 701--746.

\bibitem[KS16a]{KS2016}
S. Kawaguchi and J.~H. Silverman, Dynamical canonical heights for Jordan blocks, arithmetic degrees of orbits, and nef canonical heights on abelian varieties, Trans. Amer. Math. Soc. {\bf 368} (2016), no.~7, 5009--5035; MR3456169

\bibitem[KS16b]{KS16a}
S.~Kawaguchi and J.~H.~Silverman,
On the dynamical and arithmetic degrees of rational self-maps
of algebraic varieties,
J. Reine Angew. Math. {\bf 713} (2016), 21--48.

\bibitem[LS25]{LS2025}W. Luo and J. Song, Arithmetic degrees of dynamical systems over fields of characteristic zero, Int. Math. Res. Not. IMRN {\bf 2025}, no.~21, Paper No. rnaf247, 27 pp.; MR4979573

\bibitem[Mats57]{Mats1957}T. Matsusaka, The criteria for algebraic equivalence and the torsion group, Amer. J. Math. {\bf 79} (1957), 53--66; MR0082730

\bibitem[Mats20]{Matsuzawa2020}
Y. Matsuzawa, On upper bounds of arithmetic degrees, Amer. J. Math. {\bf 142} (2020), no.~6, 1797--1820; MR4176545

\bibitem[Mor00]{Mo00}
A. Moriwaki, Arithmetic height functions over finitely generated fields, Invent. Math. {\bf 140} (2000), no.~1, 101--142; MR1779799

\bibitem[Ohni22]{Ohnishi2022}
Tomoya Ohnishi.
Arakelov geometry over an adelic curve and dynamical systems,
Ph.D. thesis, Kyoto University, 2022.

\bibitem[Seid50]{Seidenberg1950}
A. Seidenberg, The hyperplane sections of normal varieties, Trans. Amer. Math. Soc. {\bf 69} (1950), 357--386; MR0037548

\bibitem[Song26]{Song2026}
J.~Song, A high-codimensional Yuan's inequality and its application
to higher arithmetic degrees, Adv. Math. {\bf 502} (2026),
Paper No.~111180.

\bibitem[Stacks]{stacks-project}
The Stacks Project Authors, The Stacks Project, 2026,
\url{https://stacks.math.columbia.edu}.

\bibitem[SGA~4]{SGA4}
{\it Th\'eorie des topos et cohomologie \'etale des sch\'emas. Tome 3}, Lecture Notes in Mathematics, Vol. 305, Springer-Verlag, Berlin-New York, 1973; MR0354654

\bibitem[Tru20]{Tru20}
T.~T. Truong, Relative dynamical degrees of correspondences over a field of arbitrary characteristic, J. Reine Angew. Math. {\bf 758} (2020), 139--182; MR4048444

\bibitem[Voj09]{Vojta2009}
P. Vojta, Diophantine approximation and Nevanlinna theory, in {\it Arithmetic geometry}, 111--224, Lecture Notes in Math., 2009, Springer, Berlin, ; MR2757629

\bibitem[Xie23]{Xie2023}
J.~Xie,
Remarks on algebraic dynamics in positive characteristic,
J. Reine Angew. Math. {\bf 797} (2023), 117--153.

\bibitem[XY26]{XY2026}
J. Xie and X. Yuan, Partial heights, entire curves, and the geometric Bombieri-Lang conjecture, Acta Math. {\bf 236} (2026), no.~2, 365--416; MR5104982

\end{thebibliography}
\end{document}